A Generalized Non-local Quasicontinuum Approach for Efficient Modeling of Architected Truss-based Lattice Structures

Zi Li[1], Fan Yang[2,3,4,5], Qingcheng Yang[1,5*]

[1]Shanghai Key Laboratory of Mechanics in Energy Engineering, School of Mechanics and Engineering Science, Shanghai Institute of Applied Mathematics and Mechanics, Shanghai University, Shanghai, China

[2]School of Aerospace Engineering and Applied Mechanics, Tongji University, Shanghai, China

[3]Key Laboratory of AI-Aided Airworthiness of Civil Aircraft Structures, Civil Aviation Administration of China, Tongji University, Shanghai, China

[4]State Key Laboratory of Structural Analysis, Optimization and CAE Software for Industrial Equipment, Dalian University of Technology, Dalian, China

[5]Shanghai Institute of Aircraft Mechanics and Control, Shanghai, China

\* Corresponding author, Email: qyang@shu.edu.cn, Tel: 86-18116136627.

# Abstract

To mitigate the substantial computational costs associated with modeling the mechanical behavior of large-scale architected lattice structures, this work introduces a concurrent multiscale approach: the Generalized Non-local Quasicontinuum (GNQC) method. GNQC generalizes the classical nonlocal Quasicontinuum framework by eliminating the assumption of affine or high-order deformation patterns for accurate energy sampling in coarse-grained regions and by ensuring consistency with general finite element shape functions used for coarse-graining. The introduced GNQC method offers three key features: (1) a constitutive-model-consistent framework that employs the same lattice constitutive relationship in both the locally full-resolution region and the coarse-grained domain, similar to existing nonlocal QC approaches; (2) a shape-function consistent energy sampling mechanism that aligns with the interpolation order of the generally employed shape functions, differing significantly from existing Quasicontinuum works and substantially reducing computational costs; and (3) consistent interfacial compatibility, which enables seamless energy and force transfer across interfaces between regions of different resolutions without cumbersome interfacial treatments. The performance of GNQC is validated through a series of numerical test cases—including tension, clamped bending, three-point bending, and crack propagation problems—that demonstrate good accuracy. Additionally, the error analysis and convergence behavior of GNQC are investigated.

# Introduction

With the advancement of additive manufacturing (AM) technology, the potential of architected lattice structures has been unlocked, attracting increasing attention. Lattice structures are cellular architectures characterized by non-random geometric features, formed by the periodic spatial arrangement of regular lattices[1]. This non-random characteristic endows lattice structures with highly controllable geometric properties and exceptional performance, including low elastic modulus[2,3], negative Poisson's ratio[4], large specific surface area[5], significant internal porosity[6], and tailored coefficients of thermal expansion[7]. These outstanding properties have enabled lattice structures to be utilized in a wide range of fields, such as biomedical engineering[8–10], lightweight aerospace and automotive structures[11], energy absorbers[12,13], heat exchangers[14], sound insulation materials[6,15], and catalytic applications[5,16]. Therefore, accurately and efficiently predicting the mechanical performance of lattice structures is of critical importance for their further development.

Although experiments can directly obtain the mechanical performance of lattice structures with various geometric characteristics manufactured by AM technology, experimental methods face notable limitations in engineering applications[17]. They require a large number of samples to minimize errors and improve accuracy, but the high manufacturing cost and slow production speed of AM make this try-and-error approach uneconomical.

Numerical modeling, on the other hand, offers an efficient means to predict lattice performance, facilitating faster and more cost-effective design processes while providing deeper insights into deformation mechanisms and enabling the decoupling of complex mechanical phenomena. These methods can be categorized into three groups:

homogenization methods, direct finite element analysis, and multi-scale methods. Each approach possesses distinct strengths and is particularly suited to addressing different aspects of lattice structure analysis, as detailed below.

The fundamental principle of the homogenization method is to assume scale separation to derive the effective properties of a heterogeneous material by analyzing a small, representative portion—known as the Representative Volume Element (RVE)—that encompasses the primary microstructural features. This approach circumvents the need for large-scale simulations[18–20] and is particularly well-suited for periodic lattice structures[21–29]. For example, Kumar and McDowell[23] derived a generalized continuum representation for two-dimensional periodic honeycomb solids by treating them as micropolar continua, thereby accurately capturing the overall behavior of various lattice structures. Similarly, Arabnejad and Pasini[25] investigated the mechanical properties of different lattice topologies across a wide range of relative densities, providing insights into the accuracy of homogenization schemes. Ptochos and Labeas[26] predicted the mechanical behavior of a structure comprising 1,000 unit cells using a homogenized model, significantly reducing computational costs. Moreover, Somnic and Jo[29] discussed both the potential and challenges of applying homogenization techniques to lattice materials. Although homogenized models offer high efficiency, they are limited in accurately capturing localized deformation behavior, as such deformations typically violate the scale separation assumption and require high-resolution modeling.

In contrast to the homogenization method, direct Finite Element Method (FEM) modeling enables precise prediction of the mechanical performance of lattice structures, with 3D solid elements delivering the highest accuracy[30–34]. However, because applying 3D solid elements to large-scale lattice structures is computationally expensive, beam or truss elements are often used as a more efficient alternative[35]. As early as 1982, Gibson et al.[36] employed standard beam theory to analyze the mechanical properties of two-dimensional honeycomb materials, effectively characterizing their behavior in terms of cell wall bending, elastic buckling, and plastic collapse. Later, Luxner et al.[37] modeled lattice structures using Timoshenko beam elements, achieving accurate representations of mechanical behavior while substantially reducing computational costs. Similarly, subsequent studies[38–42] have developed lattice models based on beam or truss elements, though these approaches can still entail significant computational expenses for large-scale applications.

To further reduce computational cost while accurately capturing essential local microscopic phenomena and addressing macroscopic problems, many researchers have adopted concurrent multiscale methods, which can be classified into two categories. The first approach decomposes the original system into coarse-grained domains and fine regions, preserving the discrete nature in the fine regions to explicitly capture localized deformations. In contrast, the second approach employs homogenization to establish an equivalent continuum model and then enriches it to account for local phenomena, with the resulting model solved using standard continuum techniques.

A representative example of the first category is the well-known Quasicontinuum (QC) method, originally introduced by Tadmor et al.[43] for atomistic-to-continuum coupling. In contrast, the second approach is exemplified by continuum formulations that incorporate cohesive zone models[44–46]. Mikeš et al.[47] provided a detailed description and in-depth comparison of these two techniques through three-point bending tests on specimens exhibiting crack propagation. Their study found that the locally enriched multiscale model is more efficient but lacks sufficient detail and accuracy, whereas QC-based models capture all the essential phenomena of the full model while significantly reducing the number of degrees of freedom (DOFs).

Traditional QC, also known as local QC, employs different constitutive models in the localized and coarse-grained regions. For example, the Cauchy-Born approximation[43] is used in the low-resolution domain to upscale microstructural information into a continuum constitutive model, while a more accurate discrete constitutive model is maintained in the high-resolution region. However, this dual-model approach often introduces interfacial incompatibility due to constitutive mismatches and may require cumbersome interfacial treatments.

On the other hand, nonlocal QC models eliminate interfacial incompatibility by using the same accurate discrete constitutive model in both fully resolved and coarse-grained regions, and by applying a summation rule[48–53] to efficiently compute energy and forces in the coarse-grained domain. For instance, Beex and co-workers[54–56] proposed a virtual work–based QC formulation capable of handling non-conservative lattice models. Because this approach is grounded in the principle of virtual work, its governing equations are constructed more accurately than those of force-based QC methods, where internal forces and virtual displacements of lattice nodes are not fully correlated.

The nonlocal QC framework has also been applied to reduce the computational cost of modeling plane lattices composed of Euler–Bernoulli beams[57]. However, as linear interpolation cannot accurately capture the out-of-plane deformations of these beam lattices, subsequent studies explored cubic interpolation for node displacement components and quadratic interpolation for node rotation components[58]. Comparative examples of uniaxial tension and pure bending tests demonstrated that, for an equivalent computational cost, higher-order QC methods can achieve accuracy comparable to that of linear QC methods[59].

Further advancements include the work of Rokoš et al.[60,61], who extended the variational quasi-continuum method by proposing an adaptive scheme that refines the region ahead of a crack tip while coarsening the region in its wake—thereby reducing computational time while maintaining a low relative energy error. Chen et al.[62–64]developed an adaptive coarse-graining strategy using three-dimensional co-rotational beam finite elements with embedded plastic hinges, proposing a universal refinement criterion based on the energy of specialized unit cells on the coarse-grained surface and extending the adaptive generalized quasi-continuum (AGQC) method to elastoplastic lattices. Similarly, Kochmann and co-workers[65] extended the QC method by incorporating geometric nonlinearity through a corotational beam model and by employing representative unit cells with affine interpolation for coarse-graining, enabling the simulation of nonlinear deformations in general periodic truss lattices. They further introduced a mixed-order QC method that combines linear and quadratic interpolation to address the limitations of affine interpolation in bending-dominated lattices[66]. It is worth mentioning that local QC models and most nonlocal QC approaches for architected lattice structures rely on an affine deformation assumption for accurately energy sampling in the coarse-grained region, which is only consistent with linear finite element shape functions.

This paper introduces the Generalized Non-local Quasicontinuum (GNQC) method for modeling of architected lattice structures, which eliminates assumptions of affine or high-order deformation pattern in the low-resolution region, coarse-grained by general finite element shape functions. This approach builds on the previously proposed multiresolution molecular mechanics[50,67–69]. GNQC employs the same discrete constitutive model in both the fully resolved and coarse-grained regions, as in traditional non-local QC methods. In the low-resolution region, GNQC introduces a novel energy sampling framework that is consistent with general interpolation shape function for reduction of degrees of freedom (DOFs), thereby generalizing the traditional non-local QC. The proposed energy sampling framework is akin to the role of gauss quadrature for evaluating energy integrals in conventional finite element method in continuum mechanics. Specifically, the proposed energy sampling rule first mathematically derives the energy distribution within a finite element associated with a given shape function for primary sampling, and then employs secondary sampling to refine the regions where the derived energy distribution is not followed strictly, ensuring interface compatibility. Numerical examples are employed to analyze the error and convergence behavior of GNQC as well as its performance using lattice structures under tension, bending and crack propagation conditions.

The remainder of this paper is organized as follows: Section 2 provides a comprehensive overview of the framework of GNQC. Section 3 investigates the interfacial compatibility and convergence behavior of GNQC through tensile and bending deformations employing triangular and square lattices, utilizing rigorously defined error estimates. Section 4 demonstrates the method's application to three-point bending and notched uniaxial tensile tests, highlighting its effectiveness in capturing fracture phenomena. Finally, Section 5 presents the study's conclusions and key insights, and outlines potential directions for future research.

## 2. Framework

## 2.1 Overview

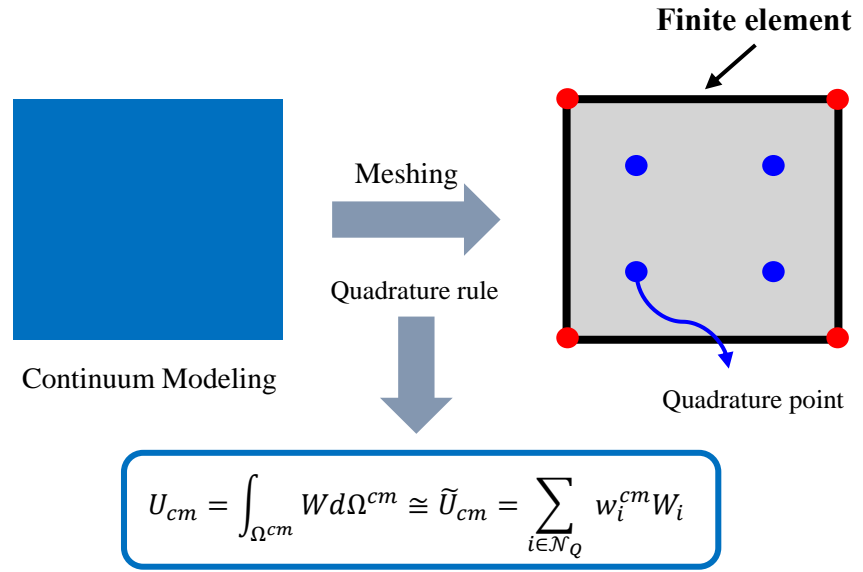

(a)

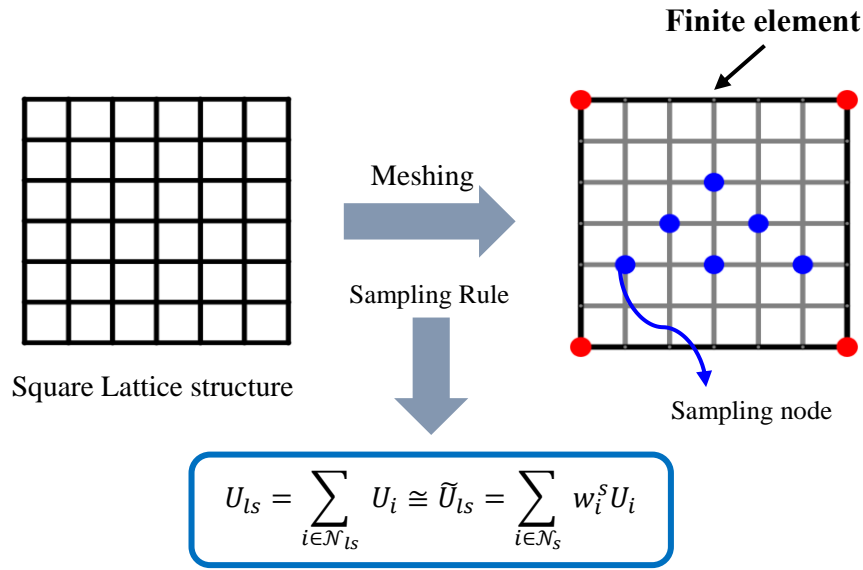

(b)

**Fig. 1.** Schematic overview of energy sampling in lattice structure modeling using GNQC by illustrating the analogy with quadrature rule in Finite Element Method:

Continuum modeling (a) and discrete lattice structure modeling (b)[69]. The red dots represent finite element nodes which are the degrees of freedom.

In this section, the overall concept of generalized nonlocal Quasicontinuum GNQC method is explained through an analogy with continuum mechanics (*cm*) and the finite element method (FEM). As shown in Fig. 1a, FEM approximates the continuum model by introducing a shape function $\phi$, enabling the total potential energy of the continuum model, $U_{cm}$ to be expressed approximately as:

$$U_{cm} = \int_{\Omega^{cm}} W d\Omega^{cm} \cong \widetilde{U}_{cm} = \sum_{i \in \mathcal{N}_Q} w_i^{cm} W_i \qquad (1)$$

where $W$ is the potential energy density, $\Omega^{cm}$ represents the continuous domain, $\mathcal{N}_Q$ refers to the index set of the quadrature points, and $w_i^{cm}$ is the corresponding weight for a quadrature point $i \in \mathcal{N}_Q$. In this study, the calligraphic letter $\mathcal{N}$ denotes the index set, while the Roman letter $N$ represents its cardinality.

To evaluate the energy integrals in Eq.(1), the Gaussian quadrature method is widely adopted. Similarly, in Fig. 1b, a finite element mesh $\phi$ is also utilized to reduce the degrees of freedom of the lattice structure (*ls*), and its potential energy $U_{ls}$ can be approximated as:

$$U_{ls} = \sum_{i \in \mathcal{N}_{ls}} U_i \cong \widetilde{U}_{ls} = \sum_{i \in \mathcal{N}_s} w_i^s U_i \qquad (2)$$

where $\mathcal{N}_{ls}$ is the index set of $N_{ls}$ nodes in the lattice structure, and $U_i$ is the node-wise energy for a node $i \in \mathcal{N}_{ls}$. $\mathcal{N}_s$ represents the index set of the selected $N_S$ sampling nodes, and $w_i^s$ is the weight associated with a sampling node $i \in N_s$. In Fig. 1, the red dots represent finite element nodes which are the degrees of freedom, while the blue dots respectively represent the quadrature points in continuum mechanics and sample nodes in the lattice structure modeling. The specific color scheme will be further explained in the next section. Additionally, any symbol with the subscript or superscript "*ls*" used in lattice structure modeling has the same physical meaning as the corresponding symbol with the subscript or superscript "*cm*" used in continuum mechanics.

In Eq.(2), it is crucial to apply an accurate energy sampling rule to efficiently evaluate the discrete energy summation in the coarse-grained region defined by a given finite element shape function $\phi$, analogous to Gaussian quadrature in continuum analysis. The essence of Gaussian quadrature lies in determining the optimal number, locations, and weights of quadrature points based on a well-established theoretical framework in continuum mechanics. Similarly, in coarse-grained discrete mechanics modeling with general finite element shape functions, an energy sampling rule must determine the optimal number, weights, and positions of sampling nodes. This aspect will be addressed in detail in Section 2.3.

## 2.2 Lattice node classification

In current multiscale modeling, the lattice structure model is decomposed into two different regions: the fully resolved region (high resolution) with localized deformations and the coarse-grained region (low resolution), as shown in Fig. 2. The lattice nodes are then classified into different groups to facilitate the introduction of the proposed energy sampling rule. From the perspective of degrees of freedoms (DOFs), the lattice nodes are classified into representative nodes (*RNs*) and ghost nodes (*GNs*). The DOFs of representative nodes are employed to interpolate those of ghost nodes through an employed shape function. Depending on whether an *RN* participates in interpolation, *RNs* are further categorized into interpolating *RNs* (*IRNs*, red) and non-interpolated *RNs* (*NIRNs*, black) in Fig. 2.

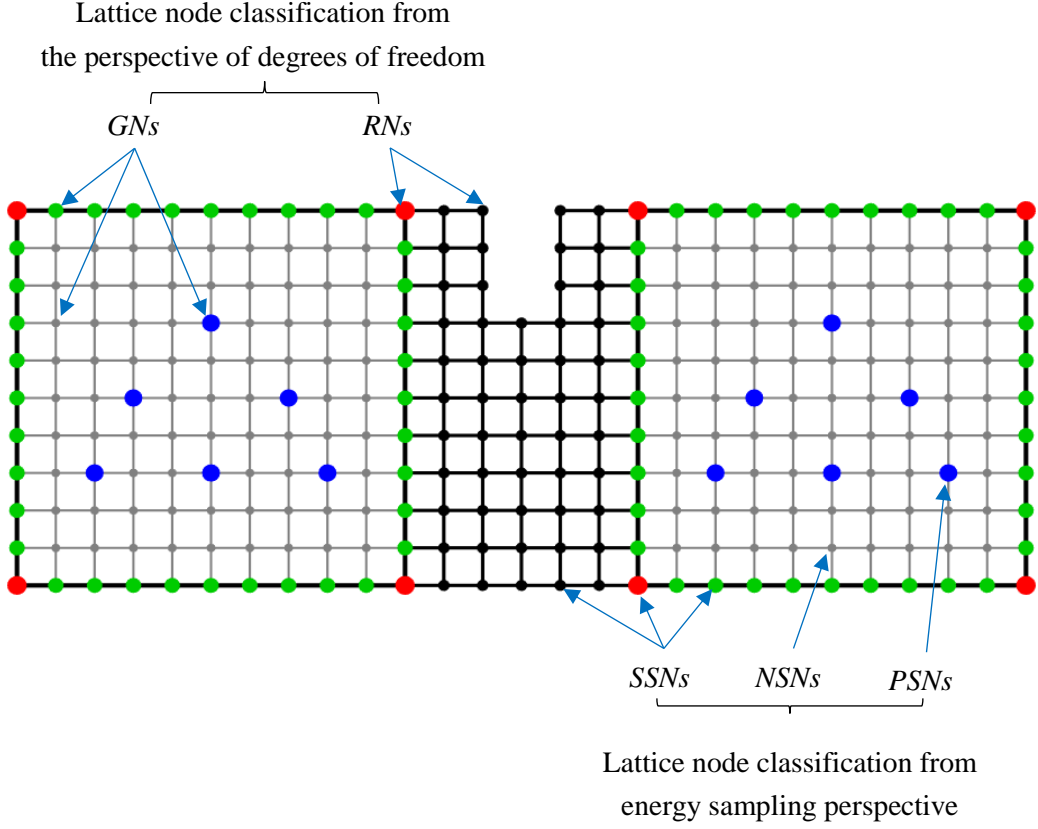

**Fig. 2.** Schematic diagram using bilinear finite element (represented by red nodes) to illustrate domain decomposition and the classification of node types in lattice structure modeling: regions occupied by black nodes represent locally full resolution domain and those occupied by finite elements denotes coarse-grained region. *RNs* and *GNs* denote representative nodes and ghost nodes, respectively. *PSNs*, *SSNs* and *NSNs* denote primary sampling nodes, secondary sampling nodes and non-sampling nodes, respectively.

From the energy sampling perspective, the lattice nodes are categorized into three types: primary sampling nodes (*PSNs*), secondary sampling nodes (*SSNs*), and non-sampling nodes (*NSNs*), as shown in Fig. 2. A *PSN i* is analogous to the Gauss quadrature point in finite element method. Its energy is used to not only sample its own but also to sample the energy of *NSNs*. Thus, the assigned weight of a *PSN* has $|w_i^s| > 1$ and that of a *NSN* has $w_i^s=0$. The *SSNs* is introduced to further improve sampling accuracy when necessary, without significantly increasing computational cost. The energy of a *SSN i* only represent its own and is explicitly calculated, rather than being sampled by *PSNs*. As such, the assigned weight of a *SSN i* has $w_i^s = 1$. It is noteworthy that representative nodes (*RNs*) belong to *SSNs* in the proposed classification.

In summary, the nodes in the lattice structure can be categorized from the perspective of degrees of freedom into *RNs* and *GNs*. On the other hand, from the perspective of energy sampling, the lattice nodes are divided into three types: *PSNs*, *SSNs*, and *NSNs*. Table 1 and Table 2 provide detailed definitions and functions for each node type. Since the lattice nodes can be classified from different perspectives, there may be overlaps among the aforementioned types. The coloring scheme used in this work to identify different node types is as follows: non-interpolating and interpolating *RNs* are represented by black and red nodes, respectively; in the case of *GNs*, primary, secondary, and non-sampling nodes are represented by blue, green, and gray nodes, respectively (Fig. 2).

**Table 1.**

Classification of lattice node types from the perspective of degrees of freedom[68].

| | | |
|---|---|---|
| **Node classification criterion:** degrees of freedom | **Rep-node:** a node whose information is the degrees of freedom of the model | **Interpolating rep-node:** rep-node that takes part in constructing finite element shape function |
| | | **Non-interpolating node:** rep-node that does not take part in interpolation |
| | **Ghost node:** a node whose information is given by interpolation | |

**Table 2.**

Classification of different node types from the perspective of energy sampling[68].

| Node classification criterion: energy sampling | **Primary sampling node (*PSN*):** a node whose energy is employed to represent the energy of itself and the non-sampling nodes | **Secondary sampling node (*SSN*):** a node whose energy is used to represent the energy of itself only and is explicitly calculated and included in the energy functional | **Non-sampling node (*NSN*):** a ghost node whose energy is represented by the energy of a primary sampling node |
|---|---|---|---|

# 2.3 Energy sampling framework

## 2.3.1 Primary sampling

In this section, we will address the questions raised in Section 2.1, namely: how to determine the optimal number, weight $w_i^s$ and locations of *PSNs* required to achieve the best accuracy and computational cost. The following derivation will also draw analogies with traditional finite element method in continuum mechanics.

It is well known that the required order of the Gaussian quadrature rule used in traditional FEM is related to the order of employed shape functions[70]. In the following derivation, we will show that the order of the proposed energy sampling rule also depends on the order of the chosen shape functions. For simplicity, we will ignore external loads and assume linear elasticity in continuum modeling. In this case, the potential energy approximation $\widetilde{U}_{cm}$ in Eq. (1) for the finite element in Fig. 1a can be expressed as follows:

$$\widetilde{U}_{cm} = \frac{1}{2}\sum_{i\in\mathcal{N}_n^{cm}} \sum_{j\in\mathcal{N}_n^{cm}} (\mathbf{u}_i^{cm})^T \mathrm{K}_{ij}^{cm}(\mathbf{u}_j^{cm}) \tag{3}$$

Here, $\mathcal{N}_n^{cm}$ refers to the index set of $N_n^{cm}$ finite element nodes, $\mathbf{u}_i^{cm}$ is the displacement vector at a finite element node $i \in \mathcal{N}_n^{cm}$, and $\mathrm{K}_{ij}^{cm}$ denotes the stiffness matrix between nodes $i$ and $j$. $\mathrm{K}_{ij}^{cm}$ can be written as:

$$\mathrm{K}_{ij}^{cm} = \int_{\Omega^{cm}} (\mathbf{B}_i^{cm})^T \mathrm{D}^{cm}(\mathbf{B}_j^{cm}) \mathrm{d}\Omega^{cm} \tag{4}$$

Here, $\mathrm{D}^{cm}$ is the material constant matrix that characterizes the constitutive behavior of the continuum, and $\mathrm{B}_i^{cm}$ is the strain matrix at the node $i$, which can be defined in a three-dimensional (3D) context as:

$$\mathrm{B}_i^{cm} = \begin{bmatrix} \phi_{i,x} & 0 & 0 \\ 0 & \phi_{i,y} & 0 \\ 0 & 0 & \phi_{i,z} \\ \phi_{i,y} & \phi_{i,x} & 0 \\ 0 & \phi_{i,z} & \phi_{i,y} \\ \phi_{i,z} & 0 & \phi_{i,x} \end{bmatrix} \tag{5}$$

From Eqs.(3)-(5), it is evident that the order of the quadrature rule needed to calculate $\widetilde{U}_c$ or $\mathrm{K}_{ij}^{cm}$ is determined by the product of the derivatives of the shape functions associated with nodes $i$ and $j$.

Similarly, for the lattice structure in Fig. 1b, it can be viewed as a truss structure with nodal interactions using truss elements. Therefore, the potential energy approximation $\widetilde{U}_{ls}$ in Eq. (2) for the finite element in Fig. 1b can be written as:

$$\widetilde{U}_{ls} = \frac{1}{2}\sum_{i\in\mathcal{N}_n^{ls}} \sum_{j\in\mathcal{N}_n^{ls}} (\mathbf{u}_i^{ls})^T \mathrm{K}_{ij}^{ls}(\mathbf{u}_j^{ls}) \tag{6}$$

where $\mathcal{N}_n^{ls}$ is the index of finite element nodes in lattice structure modeling, and the stiffness matrix $\mathrm{K}_{ij}^{ls}$ is defined as:

$$\mathrm{K}_{ij}^{ls} = \frac{1}{2}\sum_{\alpha\in\mathcal{N}_{ls}} \sum_{\beta\in\mathcal{N}_{ls}^{\alpha}} (\mathbf{B}_i^{ls})^T \mathrm{D}^{ls}(\mathbf{B}_j^{ls}) \tag{7}$$

$$\mathrm{B}_i^{ls} = [\phi_i(\mathbf{r}_{\alpha 0}) - \phi_i(\mathbf{r}_{\beta 0})]\mathbf{I} \tag{8}$$

Here, $\mathcal{N}_{ls}$ is the index set of $N_{ls}$ lattice nodes belong to the employed finite element in Fig. 1b, $\mathcal{N}_{ls}^{\alpha}$ represents the index sets of the neighboring lattice nodes interacting with the node $\alpha \in \mathcal{N}_{ls}$, and $\mathbf{r}_{\alpha 0}$ denotes the its initial position vector. $D^{ls}$ refers to the material constant or local stiffness matrix that characterizes the interaction between lattice nodes. For truss elements, $D^{ls}$=EA/$l$, where E is the Young's modulus, A is the cross-section area of the truss element, and $l$ is the truss length, which is the distance between a lattice node and its nearest neighbor. **I** is the 3 × 3 identity matrix in 3D. It is important to note that $\mathbf{B}_i^{ls}$ differs from $\mathbf{B}_i^{cm}$ because the former does not involve the derivatives of the shape functions but the difference of the shape function values evaluated at different lattice nodes.

For non-random or periodic lattice structures, the adjacent initial position vector $\mathbf{r}_{\beta 0}$ can be derived from $\mathbf{r}_{\alpha 0}$ for an interacting pair $\alpha$ and $\beta$ as follows:

$$\mathbf{r}_{\beta 0} = \mathbf{r}_{\alpha 0} + \mathbf{C} \qquad (9)$$

Here, **C** is a non-zero constant vector representing the lattice architecture and topology. By substituting Eq.

(9) into Eq.(8), we find that the evaluation of $K_{ij}^{ls}$ is governed by the product of the shape function differences computed at a lattice node $\alpha \in \mathcal{N}_{ls}$.

To better demonstrate this point, the 4-node bilinear quadrilateral element in Fig. 1b is employed as a representative example. Let $\Pi_{ij}^{ls}$ denote the product of the shape function differences as follows:

$$\Pi_{ij}^{ls} = [\phi_i(\mathbf{r}_{\alpha 0}) - \phi_i(\mathbf{r}_{\alpha 0} + \mathbf{C})][\phi_j(\mathbf{r}_{\alpha 0}) - \phi_j(\mathbf{r}_{\alpha 0} + \mathbf{C})] \qquad (10)$$

The bilinear quadrilateral shape function $\phi_i$ for a finite element node $i \in \mathcal{N}_n^{ls}$ is given as:

$$\phi_i(x, y) = a_i^0 + a_i^1 x + a_i^2 y + a_i^3 xy \qquad (11)$$

where $x$ and $y$ represent the projections of the initial position vector $\mathbf{r} = \begin{pmatrix} x \\ y \end{pmatrix}$ of an arbitrary lattice node within the element onto the *x*- and *y*- axes, respectively. The coefficient $a_i^k$ (where $k = 0-3$) are defined in terms of the finite element nodal positions used to construct $\phi_i$. Let $C_x$ and $C_y$ be the projections of the non-zero constant vector $\mathbf{C} = \begin{pmatrix} C_x \\ C_y \end{pmatrix}$ from Eq.(9) along the *x*- and *y*- directions, respectively. Then from Eq.(11), $\phi_i(\mathbf{r} + \mathbf{C})$ can be computed as:

$$\phi_i(x + C_x, y + C_y) = a_i^0 + a_i^1(x + C_x) + a_i^2(y + C_y) + a_i^3(x + C_x)(y + C_y) \quad (12)$$

The difference between $\phi_i(x, y)$ and $\phi_i(x + C_x, y + C_y)$ is given by:

$$\phi_i(x, y) - \phi_i(x + C_x, y + C_y) = b_i^0 + b_i^1 x + b_i^2 y \quad (13)$$

Here, $b_i^0$, $b_i^1$ and $b_i^2$ are coefficients that can be represented in terms of the constants $a_i^0$, $a_i^1$, $a_i^2$, $a_i^3$, $C_x$, and $C_y$. Consequently, Eq. (10) can be calculated as:

$$\Pi_{ij}^{ls} = (b_i^0 + b_i^1 x + b_i^2 y)(b_j^0 + b_j^1 x + b_j^2 y)$$

$$= c_{ij}^0 + c_{ij}^1 x + c_{ij}^2 y + c_{ij}^3 xy + c_{ij}^4 x^2 + c_{ij}^5 y^2 \quad (14)$$

Here, $c_{ij}^k$, with $k = 0-5$, is a constant defined in terms of $b_i^l$ and $b_j^l$, $l = 0-2$. Therefore, from Eqs. (6)-(8) and (14), the node-wise energy $U_\alpha$ of a lattice node $\alpha \in \mathcal{N}_{ls}$, within the corresponding bilinear element in Fig. 1b, can also be expressed as:

$$U_\alpha(\mathbf{r}_{\alpha 0}) = d_\alpha^0 + d_\alpha^1 x + d_\alpha^2 y + d_\alpha^3 xy + d_\alpha^4 x^2 + d_\alpha^5 y^2 \quad (15)$$

where $d_\alpha^k, k = 0$-5 are constants related to the finite element nodal displacement vector $\mathbf{u}_i^{ls}$, $i \in \mathcal{N}_n^{ls}$, as well as the material constants of the employed truss element and the coefficients in the corresponding bilinear shape function. Therefore, $\widetilde{U}_{ls}$ can be expressed as:

$$\widetilde{U}_{ls} = \sum_{\alpha \in \mathcal{N}_{ls}} U_\alpha \quad (16)$$

For the employed bilinear quadrilateral element in Fig. 1b, it follows from Eq. (16) that $\widetilde{U}_{ls}$ has a fully quadratic dependence on the positions of lattice nodes. As such, the functional basis to represent the potential energy function $\widetilde{U}_{ls}$ is 1, $x$, $y$, $xy$, $x^2$ and $y^2$. Thus, the optimal number of *PSN*s needed for the bilinear quadrilateral element is

6. Consequently, $U_\alpha(\mathbf{r}_{\alpha 0})$ can be expressed as in terms of the node-wise energies defined at *PSN*s, as follows:

$$U_\alpha(\mathbf{r}_{\alpha 0}) = \sum_{i \in \mathcal{N}_s} \psi_i(\mathbf{r}_{\alpha 0}) U_i, \alpha \in \mathcal{N}_{ls} \tag{17}$$

Here $\mathcal{N}_s$ is the index sets of *PSN*s, and $\psi_i$ is the Lagrange shape function or nodal shape function defined at a *PSN* $i \in \mathcal{N}_s$. It is important to note that $\psi$ is employed for energy sampling, whereas the finite element shape function $\phi$ is used for interpolating the degrees of freedoms. For the employed bilinear element, $\psi_i$ can be determined using a determinant formulation. Specifically, one defines a 6×6 Vandermonde matrix $A$ based on the nodal coordinates $(x_i, y_i)$ of *PSN*s, $i = 1,\ldots,6$, as follows:

$$A = \begin{bmatrix} 1 & x_1 & y_1 & x_1^2 & y_1^2 & x_1 y_1 \\ 1 & x_2 & y_2 & x_2^2 & y_2^2 & x_2 y_2 \\ 1 & x_3 & y_3 & x_3^2 & y_3^2 & x_3 y_3 \\ 1 & x_4 & y_4 & x_4^2 & y_4^2 & x_4 y_4 \\ 1 & x_5 & y_5 & x_5^2 & y_5^2 & x_5 y_5 \\ 1 & x_6 & y_6 & x_6^2 & y_6^2 & x_6 y_6 \end{bmatrix} \tag{18}$$

For an arbitrary lattice node $(x, y)$ within the bilinear element, the shape function $\psi_i$ associated with the *PSN* $i$ is then defined by replacing the $i$th row of $A$ with the row [1, $x, y, x^2, y^2, xy$] to form a new matrix $A_i(x, y)$. The shape function $\psi_i$ is given by

$$\psi_i(x, y) = \frac{\det A_i(x,y)}{\det A} \tag{19}$$

As seen from Eq. (19), $\psi_i$ satisfies the partition of unity and Kronecker delta properties:

$$\sum_{i \in \mathcal{N}_s} \psi_i(x, y) = 1 \tag{20}$$

$$\psi_i(x_j, y_j) = \delta_{ij}, \forall i, j \in \mathcal{N}_s \tag{21}$$

Then from Eqs. (16) and (17), the total energy $\tilde{U}_{ls}$ of the discrete lattice structure can be represented as follows:

$$\tilde{U}_{ls} = \sum_{\alpha \in \mathcal{N}_{ls}} U_\alpha = \sum_{\alpha \in \mathcal{N}_{ls}} \sum_{i \in \mathcal{N}_s} \psi_i(\mathbf{r}_{\alpha 0}) U_i = \sum_{i \in \mathcal{N}_s} \sum_{\alpha \in \mathcal{N}_{ls}} \psi_i(\mathbf{r}_{\alpha 0}) U_i \tag{22}$$

Comparing Eqs. (22) and (2), we note that:

$$w_i^s = \sum_{\alpha \in \mathcal{N}_{ls}} \psi_i(\mathbf{r}_{\alpha 0}) \tag{23}$$

It is worthy to note that for a general finite element shape function $\phi$, the process of determining the functional basis to represent $\widetilde{U}_{ls}$, identifying the optimal number of primary sampling nodes (*PSN*s), constructing the energy sampling shape function $\psi_i$ for *PSN*s, and computing the corresponding weights follows the same procedure as that applied to the employed bilinear quadrilateral shape function.

As for the locations of *PSNs*, a fundamental requirement is that determinant of the Vandermonde matrix *A* based on *PSNs* must be nonzero. Ideally, the *PSNs* should not be arranged in a way that makes Vandermonde matrix *A* ill-conditioned. For example, for the employed bilinear element, the selected *PSNs* must be non-collinear, spread out in the employed finite element. The effects of *PSNs* selection will be further discussed in Section 3 through numerical tests.

## 2.3.2 Secondary sampling

The role of secondary sampling nodes (*SSNs*) as defined in Section 2.2 is to complement the energy sampling using primary sampling nodes (*PSNs*). Note that in deriving the energy sampling rule using *PSNs* in Section 2.3.1, it is assumed that each lattice node has the same neighboring environment within the same finite element, as implicitly indicated in Eqs. (7) and (9). However, for the interpolating representative nodes (*IRNs*) acting as finite element nodes and the ghost lattice nodes near the element boundary, their neighbors may locate in different elements or even interacts with the non-interpolation nodes in the locally full-resolution region, as shown in Fig. 2. For such cases, *SSNs* are introduced to further improve energy sampling. In Section 3, different combinations of primary and secondary energy sampling schemes will be tested and recommendations will be made on the most effective energy sampling approaches to optimize both accuracy and efficiency based on numerical tests.

## 2.4 Governing equations

In this section, the governing equations of the Generalized Non-local Quasicontinuum (GNQC) method will be derived based on the variational principle, starting with the potential energy approximation $\widetilde{U}_{ls}$. Depending on the classification of the sampling nodes from the energy perspective, the potential energy $\widetilde{U}_{ls}$ can be expressed as a combination of the energies from *PSNs* and *SSNs*:

$$\widetilde{U}_{ls} = \sum_{i \in \mathcal{N}_s} w_i^s U_i = \sum_{i \in \mathcal{N}_{SSN}} U_i + \sum_{i \in \mathcal{N}_{PSN}} w_i^s U_i \tag{24}$$

where $\mathcal{N}_{SSN}$ and $\mathcal{N}_{PSN}$ are the index set of *SSNs* and *PSNs*, respectively. It is important to note that for each *SSN*, $w_i^s = 1$, and for *PSNs*, $w_i^s$ was determined using Eq. (23) in the previous section. Additionally, Eq. (23) needs to be modified to account for the introduction of *SSNs*, as follows:

$$w_i^s = \sum_{\alpha \in (\mathcal{N}_{ls} \setminus \mathcal{N}_{SSN})} \psi_i(\mathbf{r}_{\alpha 0}) \tag{25}$$

Physically, Eq. (25) means that since the *SSNs* are employed to only represent their energies explicitly, primary energy sampling should not re-sample the energies of *SSNs*. According to the variational principle, to derive the equilibrium equations, the total energy $\widetilde{U}_{ls}$ must satisfy the condition that the negative derivative with respect to the representative node displacement vector $\mathbf{u}_i^{ls}$ of a *RN* $i \in \mathcal{N}_{RN}^{ls}$, (analogous to the degrees of freedom in the finite element method), is 0. This leads to:

$$\mathbf{F}_i^{ls}(\mathbf{u}^{ls}) = -\sum_{j \in \mathcal{N}_{SSN}} \frac{\partial U_j(\mathbf{u}^{ls})}{\partial \mathbf{u}_i^{ls}} - \sum_{k \in \mathcal{N}_{PSN}} w_k^s \frac{\partial U_k(\mathbf{u}^{ls})}{\partial \mathbf{u}_i^{ls}} + \sum_{\alpha \in \mathcal{N}_{ls}} \phi_i(\mathbf{r} = \mathbf{r}_{\alpha 0}) \mathbf{f}_\alpha^{ext} = \mathbf{0}, \forall i \in \mathcal{N}_{RN}^{ls} \tag{26}$$

Here, $\mathbf{u}^{ls}$ denote the displacement vector for all the *RNs*, and $\mathbf{f}_\alpha^{ext}$ be an external force vector applied to a node $\alpha \in \mathcal{N}_{ls}$. And $\mathbf{f}_\alpha^{ext}$ is partitioned to a *RN* using the employed shape function $\phi$ as in conventional FEM. Eq.(26) is the governing equation under the proposed GNQC framework.

# 3. Interficial compatibility and convergence behavior study

In this section, to verify the interface consistency, and investigate the convergence behavior of the introduced GNQC, strict error norms and structures are first defined as in[67] Bilinear elements are used to coarse-grain the two-dimensional square and triangular lattice structures. The connections between lattice nodes are modeled using truss elements, and the material properties are consistent with those reported in the literature[71] (see Table 3).

**Table 3.**

The employed material properties using truss element.

| Young's modulus: E (MPa) | yield stress: $\sigma_y$ (MPa) | Length: $l_0$ (mm) | Cross-section area: A (mm$^2$) |
| --- | --- | --- | --- |
| 70e3 | 134 | 10 | 1 |

## 3.1 Error quantification

### 3.1.1 Sources of Error

The two types of errors that exist in traditional finite element methods are also present in our method: (1) the discretization error, which depends on the element size and the order of the interpolation shape functions, and (2) the sampling error, which is closely related to the energy sampling rule used and arises from the selection of sampling nodes or the assignment of weights[68]. As shown in Fig. 3, as a simple example, let $\mathbf{U}^{FR}$ represent the displacement solution vector obtained from the full resolution (FR) (model A in Fig. 3a), and $\mathbf{U}^{RR}$ represent the displacement solution vector obtained from the reduced resolution (RR) model (model B in Fig. 3b).

In addition, to quantify the discretization error, a special sampling model was designed (model C in Fig. 3c). In model C, after coarse-graining the full-resolution model, the energy of each node is explicitly considered, with each lattice node being treated as a secondary sampling node (*SSNs*), represented by green dots. It is important to note that representative nodes (red dots in Fig. 3) are also assigned the *SSN* node type. Since these nodes also represent degrees of freedom, they are still shown in red. In model C, only discretization error exists, with no other errors introduced. Therefore, model C is typically the most accurate model for a given finite element coarsening, though it also carries the highest computational cost after coarse-graining. Since there are no primary sampling nodes (*PSNs*) in model C, the displacement field obtained from model C is denoted as $\mathbf{U}^{FS}$, where FS indicates full sampling as all lattice nodes are treated as *SSNs*.

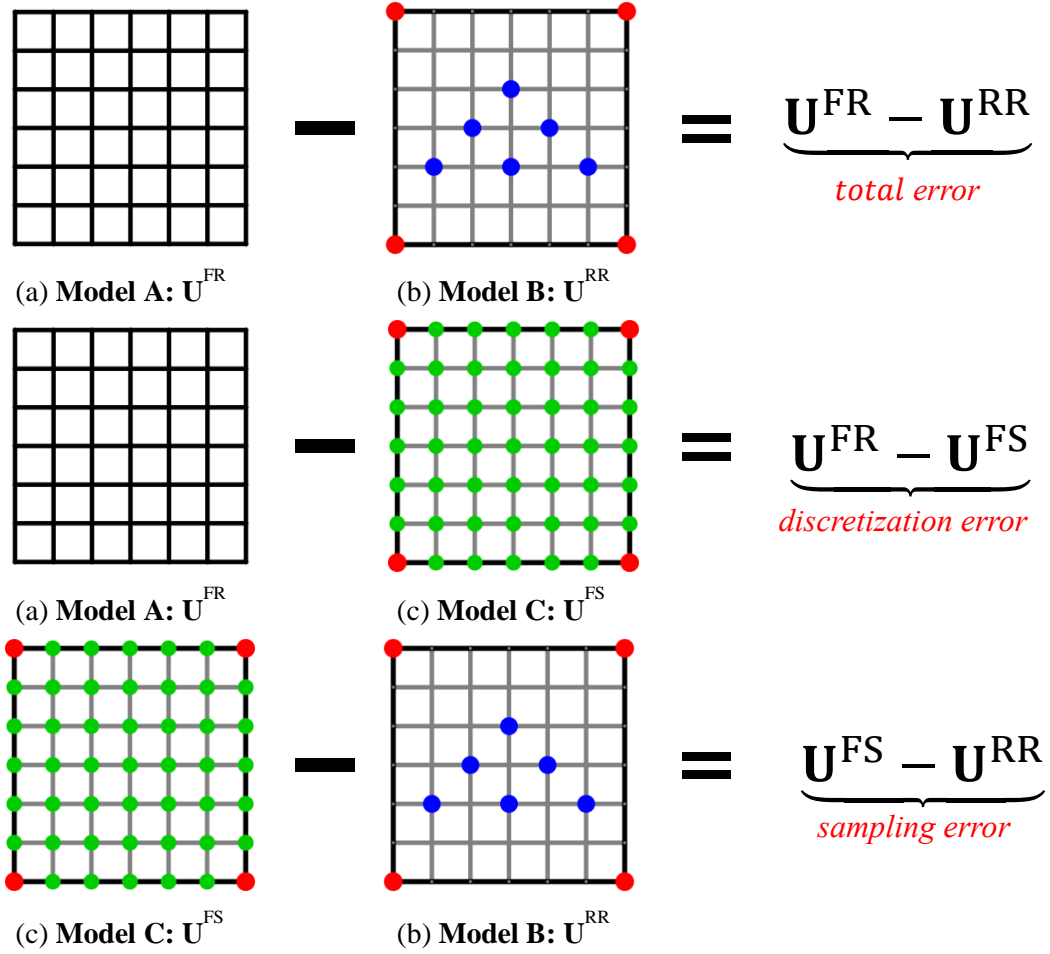

**Fig. 3.** Schematic of error sources and types: Model A: Full resolution (FR) model; Model B: Reduced resolution (RR) model with energy sampling; Model C: Full-sampling (FS) model, where the energy of each lattice node is explicitly considered. The distance between Models A and B represents the total error; the distance between Models A and C represents the discretization error; the distance between Models C and B represents the sampling error. The discretization error is fixed for a given mesh[69].

Therefore, for the displacement field, the distance between Model A and Model B can be decomposed into discretization error and sampling error, as follows:

$$\underbrace{\mathbf{U}^{FR} - \mathbf{U}^{RR}}_{\text{total error}} = \underbrace{\mathbf{U}^{FR} - \mathbf{U}^{FS}}_{\text{discretization error}} + \underbrace{\mathbf{U}^{FS} - \mathbf{U}^{RR}}_{\text{sampling error}} \quad (27)$$

## 3.1.1 Error quantification

In this work, the error is quantified using the error norms defined in reference[68]. These error norms are defined in a way that is consistent with the corresponding error norms used in continuum methods, namely the widely used $L_2$ norm and $H_1$ norm, which are used to measure the displacement and energy errors, respectively. The definitions of the displacement error norm and energy error norm are given below:

$$\| e_{disp} \| = \left[ \frac{\sum_i^{N_{ls}} (\mathbf{u}_i^{FR}-\mathbf{u}_i^{RR})^T (\mathbf{u}_i^{FR}-\mathbf{u}_i^{RR})}{\sum_i^{N_{ls}} (\mathbf{u}_i^{FR})^T \mathbf{u}_i^{FR}} \right]^{1/2} = \frac{\|\mathbf{U}^{FR}-\mathbf{U}^{RR}\|_2}{\|\mathbf{U}^{FR}\|_2} \qquad (28)$$

$$\| e_U \| = \left[ \frac{\sum_i^{N_{ls}} \sum_j^{n_i} \left((\mathbf{r}_{ij}^{FR}-\mathbf{r}_{ij}^0)-(\mathbf{r}_{ij}^{RR}-\mathbf{r}_{ij}^0)\right)^T \left((\mathbf{r}_{ij}^{FR}-\mathbf{r}_{ij}^0)-(\mathbf{r}_{ij}^{RR}-\mathbf{r}_{ij}^0)\right)}{\sum_i^{N_{ls}} \sum_j^{n_i} \left((\mathbf{r}_{ij}^{FR}-\mathbf{r}_{ij}^0)\right)^T \left((\mathbf{r}_{ij}^{FR}-\mathbf{r}_{ij}^0)\right)} \right]^{1/2} \qquad (29)$$

where $\mathbf{u}_i^{FR}$ and $\mathbf{u}_i^{RR}$ are the displacement vectors of lattice node $i$ from the full resolution model and from the different reduced-resolution models, respectively; $n_i$ is the number of interacting neighbors of a lattice node $i$. $\mathbf{r}_{ij}^{FR}$ and $\mathbf{r}_{ij}^{RR}$ represent the distance vectors between each neighboring pair, and $\mathbf{r}_{ij}^0$ represents its initial distance vector. It is important to note that all errors are normalized using the solution from the full resolution model.

Then, using the error norms defined in Eq. (28), the discretization error $e_{disp}^{disc}$ and sampling error $e_{disp}^{sam}$ defined in Eq. (27) for the displacement field can be derived:

$$\|e_{disp}^{disc}\| = \left[ \frac{\sum_i^{N_{ls}} (\mathbf{u}_i^{FR}-\mathbf{u}_i^{FS})^T (\mathbf{u}_i^{FR}-\mathbf{u}_i^{FS})}{\sum_i^{N_{ls}} (\mathbf{u}_i^{FR})^T \mathbf{u}_i^{FR}} \right]^{1/2} = \frac{\|\mathbf{U}^{FR}-\mathbf{U}^{FS}\|_2}{\|\mathbf{U}^{FR}\|_2} \qquad (30)$$

$$\|e_{disp}^{sam}\| = \left[ \frac{\sum_i^{N_{ls}} (\mathbf{u}_i^{FS}-\mathbf{u}_i^{RR})^T (\mathbf{u}_i^{FS}-\mathbf{u}_i^{RR})}{\sum_i^{N_{ls}} (\mathbf{u}_i^{FR})^T \mathbf{u}_i^{FR}} \right]^{1/2} = \frac{\|\mathbf{U}^{FS}-\mathbf{U}^{RR}\|_2}{\|\mathbf{U}^{FR}\|_2} \qquad (31)$$

Similarly, the discretization error $e_U^{disc}$ and the sampling error $e_U^{sam}$ in energy can be defined as:

$$\|e_U^{disc}\| = \left[\frac{\sum_i^{N_{ls}} \sum_j^{n_i} \left((\mathbf{r}_{ij}^{FR}-\mathbf{r}_{ij}^0)-(\mathbf{r}^{FS}-\mathbf{r}_{ij}^0)\right)^T \left((\mathbf{r}_{ij}^{FR}-\mathbf{r}_{ij}^0)-(\mathbf{r}^{FS}-\mathbf{r}_{ij}^0)\right)}{\sum_i^{N_{ls}} \sum_j^{n_i} \left((\mathbf{r}_{ij}^{FR}-\mathbf{r}_{ij}^0)\right)^T \left((\mathbf{r}_{ij}^{FR}-\mathbf{r}_{ij}^0)\right)}\right]^{1/2} \quad (32)$$

$$\|e_U^{sam}\| = \left[\frac{\sum_i^{N_{ls}} \sum_j^{n_i} \left((\mathbf{r}^{FS}-\mathbf{r}_{ij}^0)-(\mathbf{r}_{ij}^{RR}-\mathbf{r}_{ij}^0)\right)^T \left((\mathbf{r}^{FS}-\mathbf{r}_{ij}^0)-(\mathbf{r}_{ij}^{RR}-\mathbf{r}_{ij}^0)\right)}{\sum_i^{N_{ls}} \sum_j^{n_i} \left((\mathbf{r}_{ij}^{FR}-\mathbf{r}_{ij}^0)\right)^T \left((\mathbf{r}_{ij}^{FR}-\mathbf{r}_{ij}^0)\right)}\right]^{1/2} \quad (33)$$

The following relationship between the errors can be derived from the triangle inequality:

$$e_{disp} \leq e_{disp}^{disc} + e_{disp}^{sam} \quad (34)$$

$$e_U \leq e_U^{disc} + e_U^{sam} \quad (35)$$

Eqs. (34) and (35) give the upper bound of the total errors in displacement and energy error norms, respectively. The errors defined in this section will be used to quantify the accuracy of different energy sampling schemes in GNQC in different numerical tests.

## 3.2 Interfacial compatibility

In this section, a square lattice with size of $60l_0 \times 60l_0$ is used, with $l_0 = $ 10mm representing the characteristic truss length. For boundary conditions, lattice nodes at the bottom are fully fixed in both the *x* and *y* directions, while a displacement load is applied in the positive *y*-direction at the top. To test interfacial consistency, a locally full resolution (FR) model is reserved (black region in Fig. 4), and the rest is coarsened using the 4-node bilinear elements with an element size of $20l_0$. Exploiting the flexibility of secondary energy sampling, in addition to the full sampling approach (Fig. 4a), three different combinations of primary sampling (blue dots in Fig. 4) and secondary samplings (green dots in Fig. 4) are designed and compared against the full resolution model. Specifically, Scheme 1(Fig. 4b) treats all boundary lattice nodes of finite elements as secondary sampling nodes, and is termed Edge Secondary Sampling(ESS); Scheme 2 (Fig. 4c) only applies secondary sampling at the interface between coarse and full resolution domains, and is called interface secondary sampling (ISS), and Scheme 3 (Fig. 4d) employs no secondary sampling (NSS). Note that the *PSN*s are selected based on the derivations in Section 2.

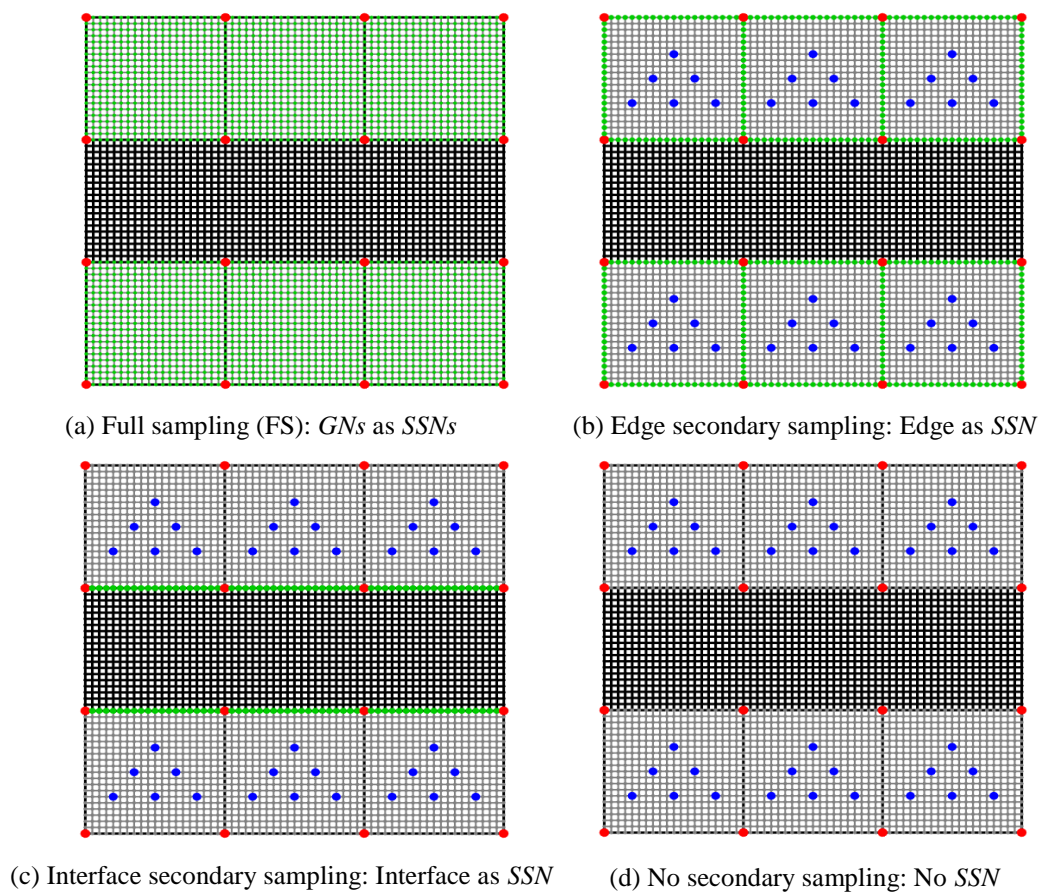

(a) Full sampling (FS): *GNs* as *SSNs*    (b) Edge secondary sampling: Edge as *SSN*

(c) Interface secondary sampling: Interface as *SSN*    (d) No secondary sampling: No *SSN*

**Fig. 4.** Schematic of different energy sampling schemes for the square lattice coarse-grained by bilinear elements: (a) Full sampling (FR): Ghost nodes (*GNs)* as *SSNs*; (b) Edge secondary sampling (ESS): Edge as *SSN*, (c) Interface secondary sampling (ISS): Interface as *SSN*, and (d) No secondary sampling (NSS): No *SSN*.

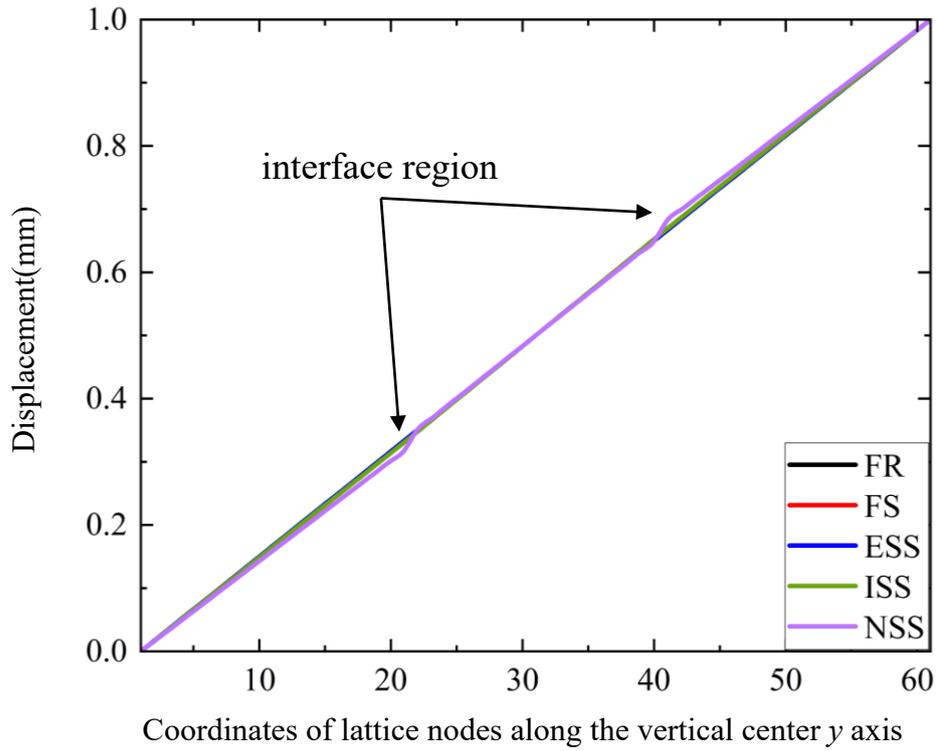

**Fig. 5.** Displacement field along the vertical center *y*-axis from the full resolution (FR) model, full-sampling (FS) model, and GNQC with edge secondary sampling(ESS), interface secondary sampling(ISS) and no secondary sampling(NSS).

Fig. 5 shows the displacement distributions across the vertical center *y*-axis for GNQC equipped different energy sampling schemes, as well as for the full resolution (FR) model and the full-sampling (FS) model, respectively. The displacement fields from the respective FS, ESS, and ISS closely match those from the FR model, indicating that the interfacial compatibility is satisfied and linear consistency is reproduced. In contrast, the displacement field from no secondary sampling (NSS) exhibits some fluctuations near the interface, which result from the sharp resolution change between the locally fully-resolved region and the coarse-grained domain. Nevertheless, these fluctuations are confined to the interface region, with small error magnitudes, and the overall accuracy is still satisfactory.

Furthermore, Table 4 quantifies various errors for the overall model under different sampling schemes, including relative displacement errors ($e_{disp}$, $e_{disp}^{sam}$, $e_{disp}^{disc}$) and relative energy errors ($e_U$, $e_U^{sam}$, $e_U^{disc}$). As mentioned in the previous section, for a given coarse-graining, the discretization error is fixed, while the sampling error is the key metric for evaluating the performance of different energy sampling rules, with FS, ESS, ISS, and NSS having a reducing order of accuracy and an increasing order of computational cost. It can be seen that the errors for ESS and ISS remain at a very low level with sampling errors less than 0.5%. In contrast, NSS exhibits relatively larger sampling errors.

**Table 4.**

Errors in displacement and energy fields for the 2D square lattice stretching example using bilinear quadrilateral elements.

| Sampling approach | Relative displacement errors | | | Relative energy errors | | |
|---|---|---|---|---|---|---|
| | $e_{disp}$ | $e_{disp}^{sam}$ | $e_{disp}^{disc}$ | $e_U$ | $e_U^{sam}$ | $e_U^{disc}$ |
| Full sampling (FS) | 2.16e-7 | 0 | 2.16e-7 | 3.42e-7 | 0 | 3.42e-7 |
| Edge secondary sampling (ESS) | 2.16e-7 | 5.94e-15 | 2.16e-7 | 3.42e-7 | 8.34e-15 | 3.42e-7 |
| Interface secondary sampling (ISS) | 0.25% | 0.25% | 2.16e-7 | 0.35% | 0.35% | 3.42e-7 |
| No secondary sampling (NSS) | 1.29% | 1.29% | 2.16e-7 | 1.86% | 1.86% | 3.42e-7 |

Fig. 6 shows a comparison of the full displacement field distributions in *y*-direction from the respective full resolution (FR) model and the GNQC with ESS, ISS, and NSS sampling schemes. The results indicate that the employed sampling schemes can reproduce the full resolution results both qualitatively and quantitatively, thereby demonstrating the flexibility of GNQC and ensuring consistency at the interface between the coarse-grained and fully-resolved domains. Based on the observations, it is recommended that interface secondary sampling (ISS) be employed in the proposed GNQC framework to achieve the optimized balance between accuracy and computational cost in practical applications.

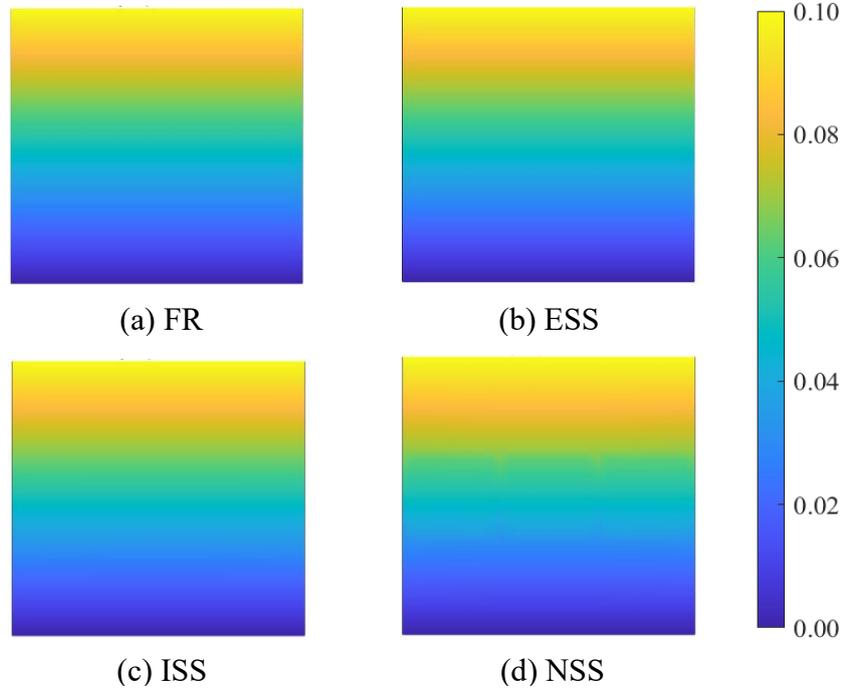

**Fig. 6.** Displacement distribution in the *y*-direction from full resolution (FR) model (a), edge secondary sampling (ESS)(b), interface secondary sampling (ISS)(c), and no secondary sampling (NSS)(d).

As mentioned in Section 2.3.1, a fundamental requirement for positioning primary sampling nodes (*PSNs*) is to ensure that their distribution must avoid rendering the Vandermonde matrix A ill-conditioned. Once this requirement is met, then the specific selection of *PSNs* should not significantly affect the accuracy, according to the proposed sampling framework in Section 2.3.1. To validate this assertion, two different selections of *PSNs* combined with the recommended interface secondary sampling (ISS) are compared, as shown in Fig. 7. Table 5 quantifies the errors of the two different selections, which yield almost identical results, with discrepancies only on the order of $10^{-13}$. This outcome supports the robustness of the proposed sampling framework. It is further recommended that *PSNs* are distributed thought out the bulk region of the employed finite element such that they are not coplanar[68].

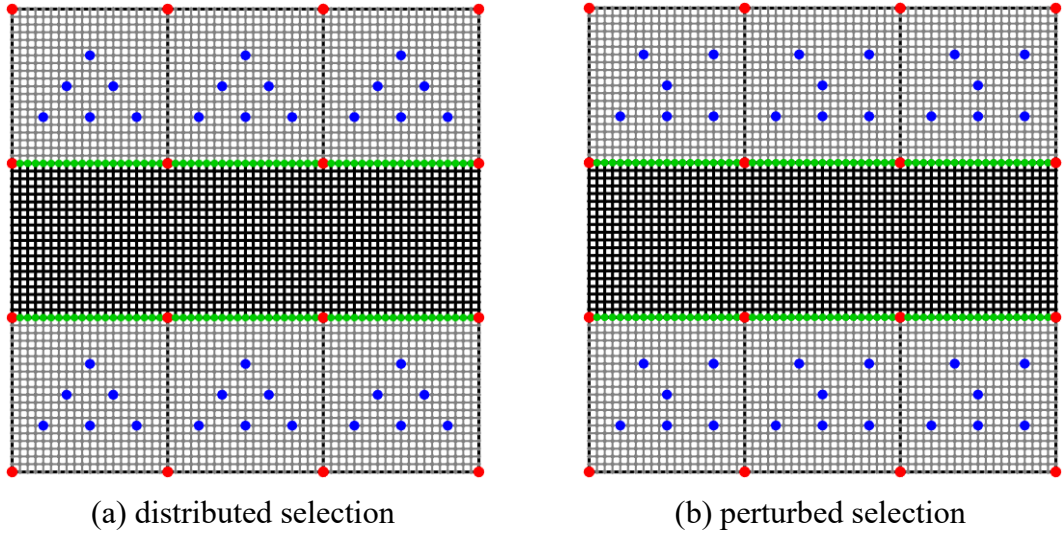

(a) distributed selection        (b) perturbed selection

**Fig. 7.** Two different selections of primary sampling nodes equipped with interface secondary sampling scheme in the proposed (GNQC).

**Table 5.**

Errors in the displacement and energy fields using GNQC with interface secondary sampling for two different selections of primary sampling node.

| Sampling rules | Relative displacement errors | | | Relative energy errors | | |
|---|---|---|---|---|---|---|
| | $e_{disp}$ | $e_{disp}^{sam}$ | $e_{disp}^{disc}$ | $e_U$ | $e_U^{sam}$ | $e_U^{disc}$ |
| ISS with distributed selection of *PSNs* | 0.25% | 0.25% | 2.16e-7 | 0.35% | 0.35% | 3.42e-7 |
| ISS with perturbed selection of *PSNs* | 0.25% | 0.25% | 2.16e-7 | 0.35% | 0.35% | 3.42e-7 |

# 3.3 Convergence study: triangular lattice under tension

In the introduced GNQC, one benefit of employing the discrete constitutive model in the low resolution region is that it automatically converges to the full resolution model when the coarse-graining mesh researches the lattice resolution. To test the convergence behavior of GNQC, a two-dimensional triangular lattice structure of size $48l_0 \times 24\sqrt{3}l_0$ (containing 2377 lattice nodes with a lattice spacing $l_0$) was subjected to a tensile load. In this setup, the bottom nodes were fixed in both the *x* and *y* directions, while a displacement load was applied in the positive *y*-direction to the top nodes. The lattice structure was coarse-grained using four bilinear elements, with an element size of $24l_0$, as shown in Fig. 8. Since there is no interface between the coarse-grained and fully-resolved regions, an additional energy sampling scheme is introduced. This scheme employs neighbors of representative nodes (*RNs*) as secondary sampling node, and is referred neighbor-assisted sampling (NAS), as shown in Fig. 8c. Compared to FS and ESS, NAS reduces the number of sampling nodes while improving the energy sampling accuracy near the interpolation nodes. The design of all employed sampling schemes is depicted in Fig. 8.

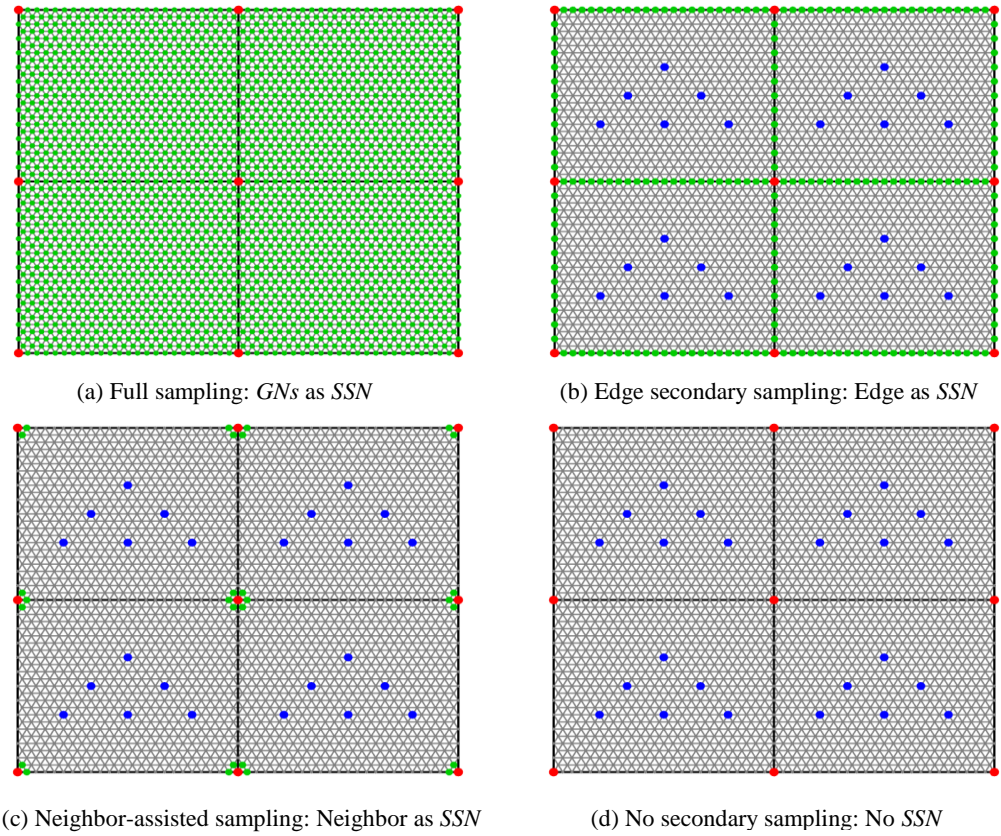

(a) Full sampling: *GNs* as *SSN*  (b) Edge secondary sampling: Edge as *SSN*

(c) Neighbor-assisted sampling: Neighbor as *SSN*  (d) No secondary sampling: No *SSN*

**Fig. 8.** Schematic of different sampling schemes for the triangular lattice coarse-grained by bilinear elements. (a) Full sampling (FR): Ghost nodes (*GNs*) as *SSNs*; (b) Edge secondary sampling (ESS): Edge as *SSN*, (c) Neighbor-assisted sampling (NAS): Neighbor as *SSN*, and (d) No secondary sampling (NSS): No *SSN*.

The *x*-displacement and *y*-displacement fields obtained from the full resolution model and GNQC with ESS, NAS and NSS are shown in Fig. 9 and Fig. 10, respectively. It is evident that the displacement field from all sampling schemes closely match those from the FR model. This demonstrates the consistency of the proposed sampling schemes in capturing the distribution characteristics of the displacement field for different lattice structures.

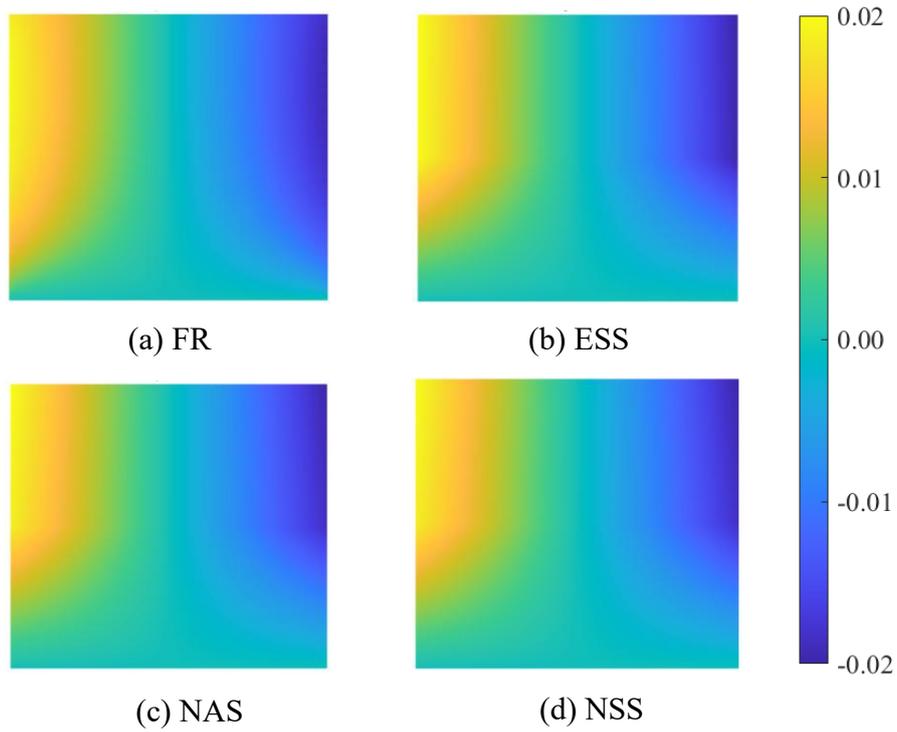

**Fig. 9.** Comparison of *x*-displacement filed from the respective full resolution (FR) model (a), edge secondary sampling (ESS) model (b), neighbor-assisted model (NAS) model (c) and no secondary sampling (NSS) model (d).

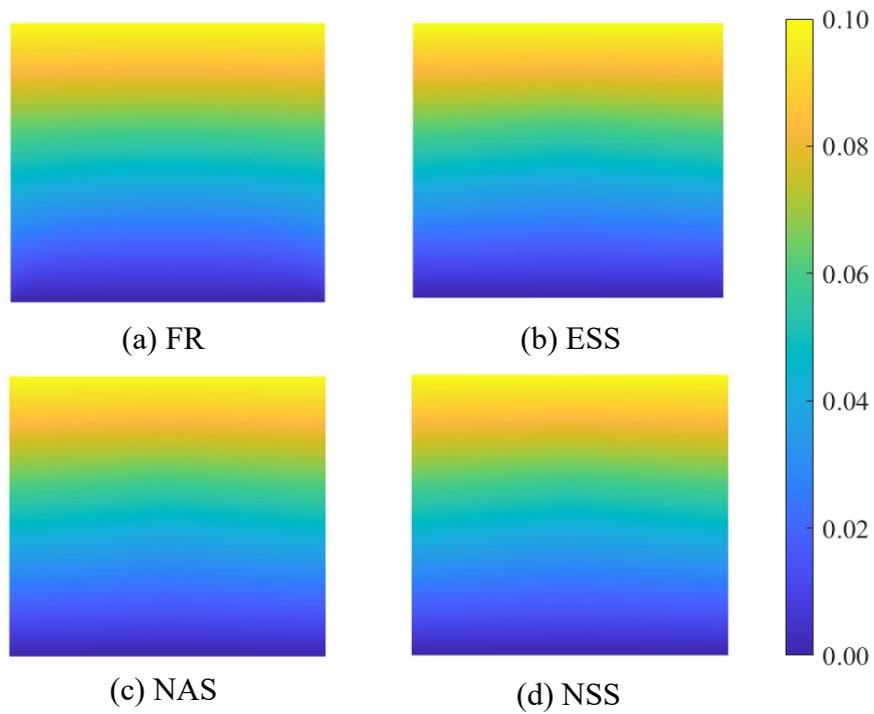

**Fig. 10.** Comparison of *y*-displacement filed from the respective full resolution (FR)

model (a), edge secondary sampling (ESS) model (b), neighbor-assisted model (NAS) model (c) and no secondary sampling (NSS) model (d).

Additionally, Table 6 quantified the relative displacement errors and energy errors for different sampling schemes. It is seen that the major contribution to the total error comes from discretization errors, indicating that the discretization error dominates sampling error.

**Table 6.**

Errors of different sampling schemes in displacement and energy fields for the 2D triangular lattice under tension using bilinear elements.

| Sampling schemes | Relative displacement errors | | | Relative energy errors | | |
|---|---|---|---|---|---|---|
| | $e_{disp}$ | $e_{disp}^{sam}$ | $e_{disp}^{disc}$ | $e_U$ | $e_U^{sam}$ | $e_U^{disc}$ |
| Full sampling (FS) | 2.35% | 0 | 2.35% | 1.12% | 0 | 1.12% |
| Edge secondary sampling (ESS) | 2.37% | 0.14% | 2.35% | 1.11% | 0.06% | 1.12% |
| Neighbor-assisted sampling (NAS) | 2.73% | 0.95% | 2.35% | 2.06% | 1.18% | 1.12% |
| No secondary sampling (NSS) | 2.78% | 1.03% | 2.35% | 2.15% | 1.28% | 1.12% |

To numerically assess the convergence behavior of GNQC, varying element sizes (see Fig. 11) were used to coarse-grain the lattice structure under tension. Fig. 12 presents the error measurements for each sampling scheme as a function of element size, using the $L_2$ and $H_1$ norms defined in Section 3.1.1. The results indicate that all sampling schemes exhibit a monotonic convergence behavior. Notably, the NSS scheme offers the highest computational efficiency but delivers the relatively lowest accuracy, while the FS scheme achieves the highest accuracy at a substantially higher computational cost. The NAS scheme strikes a balance between these extremes, and the ESS scheme attains nearly the same accuracy as the FS scheme but with a significantly reduced computational cost.

This behavior can be attributed to the number of secondary sampling nodes employed: FS utilizes the greatest number of sampling nodes, leading to high accuracy at the expense of efficiency, whereas NSS uses no secondary sampling nodes, resulting in a relatively lower accuracy but much enhanced efficiency. In addition, ESS employs fewer sampling nodes than FS yet more than NAS, which explains its intermediate trade-off between accuracy and computational cost. The measured converge orders in $L_2$ norm are 1.4011 ($R^2$ = 0.9997) for FS and 1.3857 ($R^2$=0.9998) for ESS, while in the $H_1$ norm, they are 1.3659 ($R^2$= 0.9990) and 1.3870 ($R^2$=0.9987), respectively. It should be noted that the obtained solutions are all based on discrete mechanics simulation; hence, their solutions are not as smooth as those derived from continuum models, and a conventional FEM-like convergence behavior is not observed.

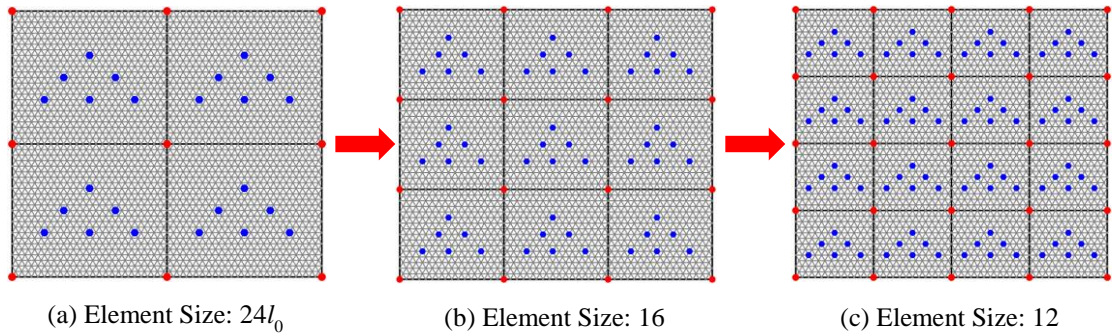

(a) Element Size: $24l_0$    (b) Element Size: 16    (c) Element Size: 12

**Fig. 11.** Diagram of triangular lattice structure coarse-grained with element sizes (a)$24l_0$, (b)$16l_0$, and (c)$12l_0$, respectively.

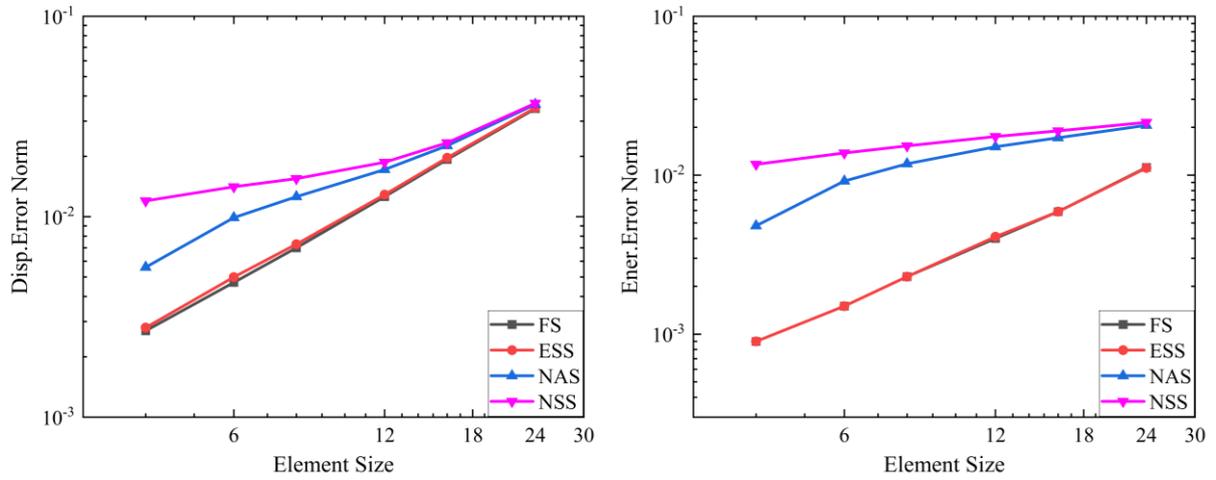

(a) relative displacement error norm  (b) relative energy error norm

**Fig. 12.** Convergence behavior of the proposed GNQC equipped with different sampling schemes using the 2D triangular lattice model under tensile: (a) relative displacement error norm and (b) relative energy error norm.

# 3.4 Convergence study: triangular lattices under bending

To further investigate the convergence behavior of GNQC, a two-dimensional lattice beam under bending is analyzed. The beam is modeled as a triangular lattice with dimensions $96l_0 \times 12\sqrt{3}l_0$, comprising a total of 2413 nodes. As for boundary conditions, the left end is fixed, while a displacement load in the negative *y*-direction is applied at the right end (Fig. 13a). The lattice structure was again coarse-grained using four bilinear elements, with an element size of $24l_0$, and the employed sampling schemes (ESS, NAS and NSS) are shown in Fig. 13b, c and d.

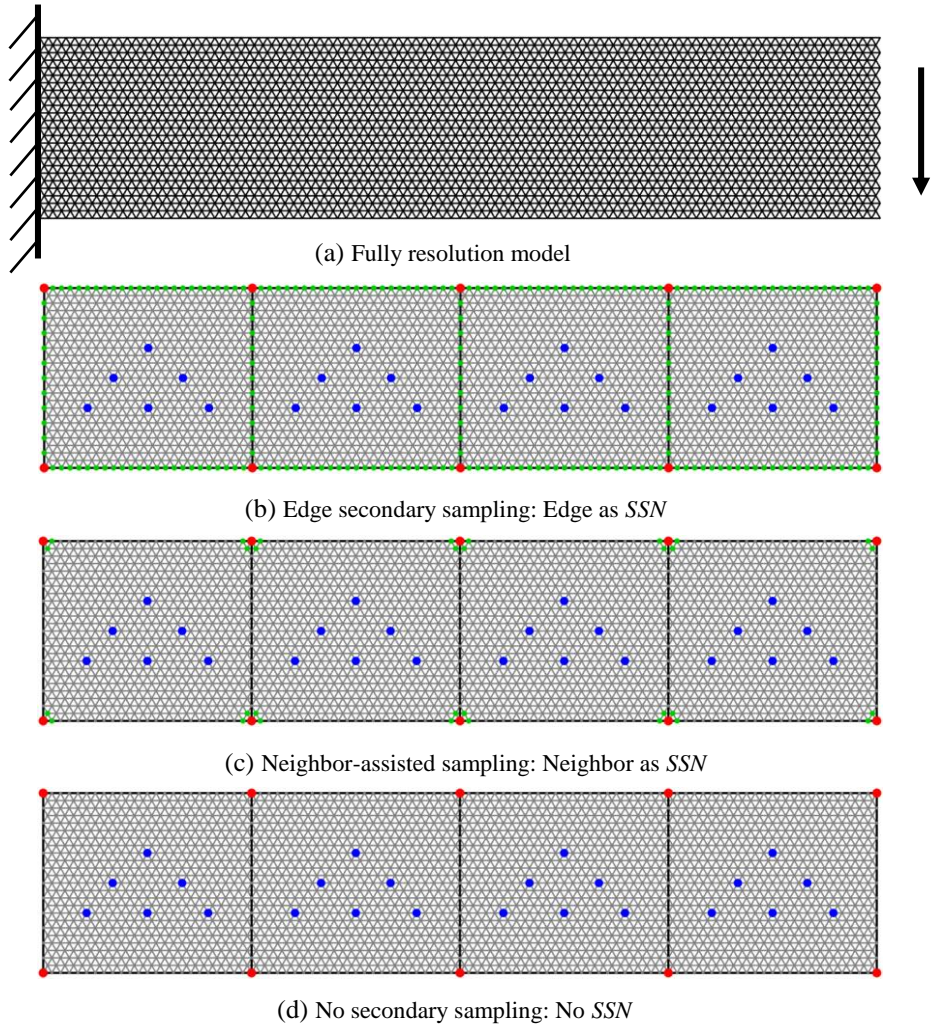

**Fig. 13.** Schematic of the triangular lattice beam under bending: (a) problem setup, and sampling schemes for the beam coarse-grained using bilinear elements: (b) Edge Secondary Sampling (ESS), (c) Neighbor-Assisted Sampling (NAS), and (d) No Secondary Sampling (NSS).

Fig. 14 compares the vertical displacement fields obtained from each sampling scheme with those from the full resolution calculation. The results indicate that each sampling scheme accurately captures the displacement distribution. Table 7 quantifies the errors for different sampling schemes, showing the discretization error dominates over the sampling error, thereby demonstrating the effectiveness of the proposed sampling framework in GNQC.

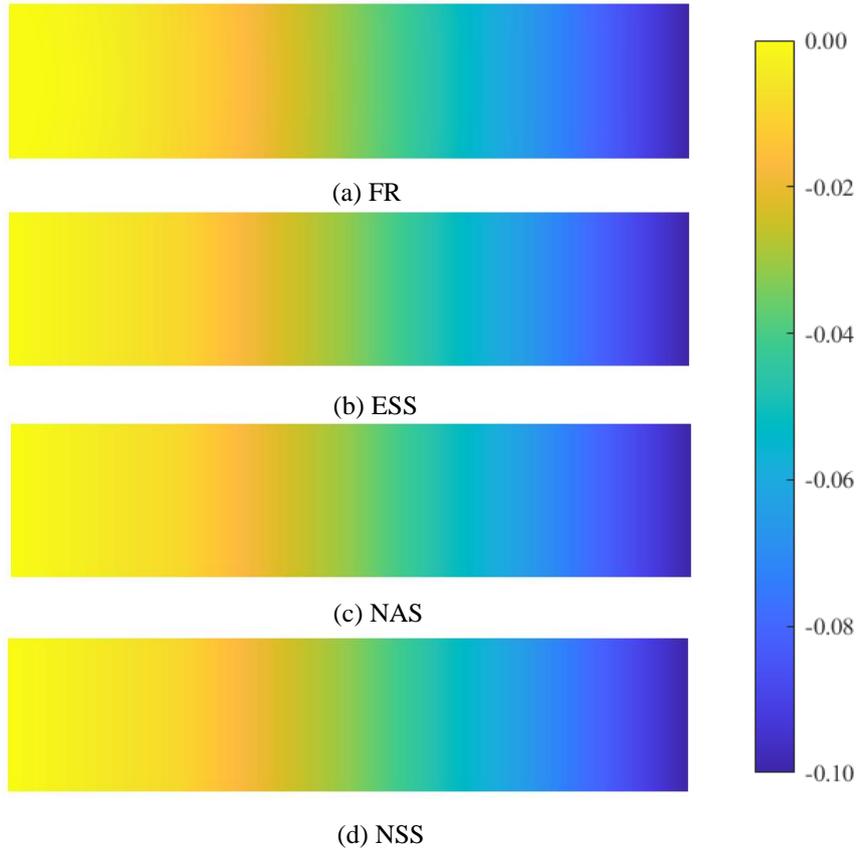

**Fig. 14.** Comparison of vertical displacement distributions in 2D beam bending from full resolution (FR) (a), edge secondary sampling (ESS) (b), neighbor assisted sampling (NAS )(c) and no secondary sampling (NSS) (d).

**Table 7.**
Errors in displacement and energy fields for the 2D triangular lattice bending example using bilinear quadrilateral elements.

| Sampling schemes | Relative displacement errors | | | Relative energy errors | | |
| --- | --- | --- | --- | --- | --- | --- |
| | $e_{disp}$ | $e_{disp}^{sam}$ | $e_{disp}^{disc}$ | $e_U$ | $e_U^{sam}$ | $e_U^{disc}$ |
| Full sampling (FS) | 1.74% | 0 | 1.74% | 0.89% | 0 | 0.89% |
| Edge secondary sampling (ESS) | 1.70% | 0.05% | 1.74% | 0.80% | 0.06% | 0.89% |
| Neighbor-assisted sampling (NAS) | 1.87% | 0.29% | 1.74% | 0.96% | 0.29% | 0.89% |
| No secondary sampling (NSS) | 2.01% | 0.31% | 1.74% | 1.03% | 0.31% | 0.89% |

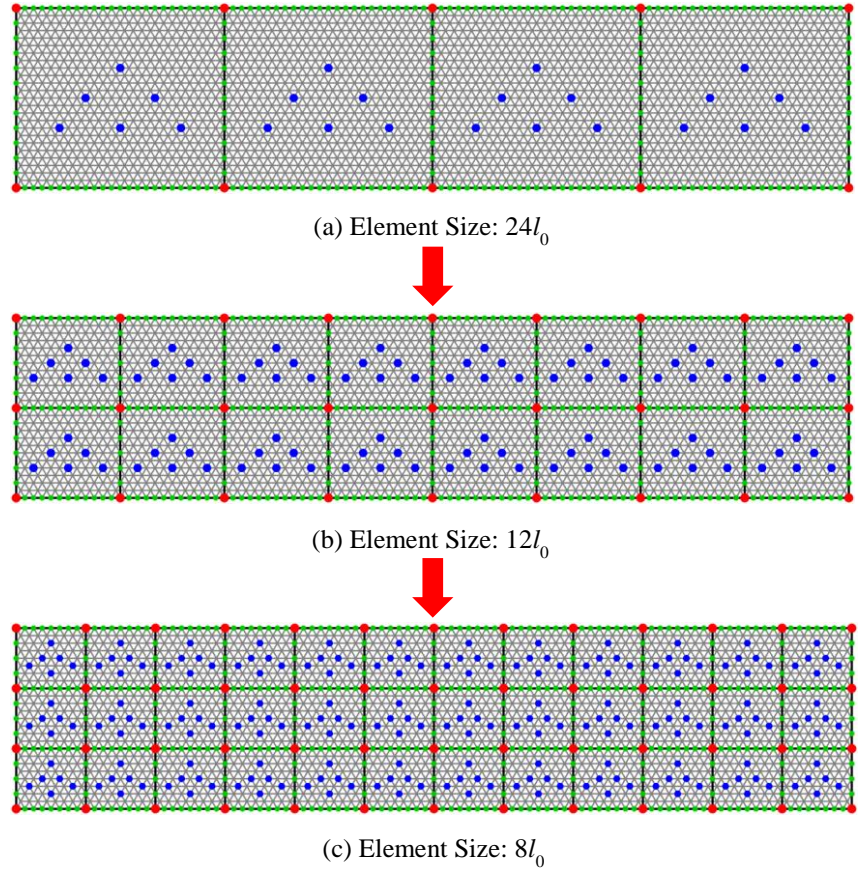

(a) Element Size: $24l_0$

(b) Element Size: $12l_0$

(c) Element Size: $8l_0$

**Fig. 15.** Diagram of triangular lattice structure coarse-grained with element sizes (a)$24l_0$, (b)$12l_0$, and (c)$8l_0$, respectively.

As shown in Fig. 15, GNQC with different element sizes is employed to investigate its convergence behavior under bending. Fig. 16 shows the convergence treads for each sampling scheme in terms of the defined $L_2$ and $H_1$ norms. Similar trends to those observed in Fig. 12 for the tension case are evident, reaffirming the consistent performance of GNQC under both tension and bending. Notably, when the element size is small, most lattice node neighbors are located in different finite elements, as shown in Fig. 15. Under such conditions, the number of secondary sampling nodes may exceed that of primary sampling nodes, eliminating the benefit of primary sampling. Therefore, to balance accuracy and efficiency, it is recommended that the element size remain relatively large. For example, the characteristic element size should be approximately 10 times the characteristic length $l_0$ of the lattice structure or larger.

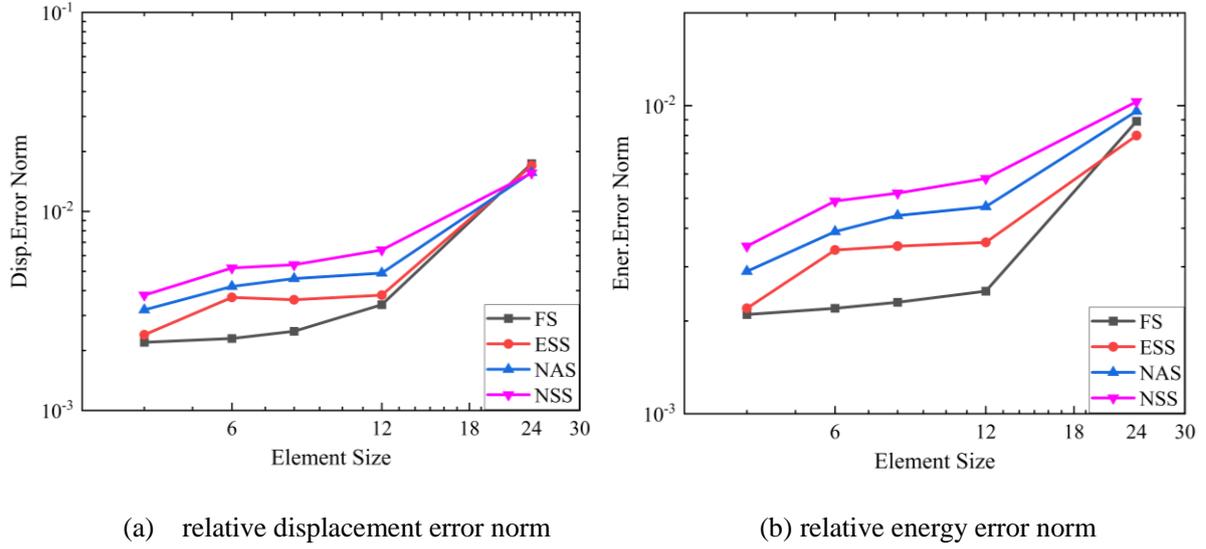

(a) relative displacement error norm　　　(b) relative energy error norm

**Fig. 16.** Convergence analysis of GNQC with different sampling schemes using a 2D triangular lattice beam under bending: (a) relative displacement error norm and (b) relative energy error norm.

# 4. Crack propagation

In this section, two numerical examples are presented to evaluate the accuracy of the proposed GNQC in capturing crack propagation in architected lattice structures. Since truss elements are used for the interactions between lattice nodes, the maximum principal stress criterion is employed to assess strut failure. The material properties remain consistent with those in Section 3, with a yield stress of 134 MPa.

## 4.1 Three-Point Bending

The mechanical behavior of a two-dimensional triangular lattice model with a notch under three-point bending is simulated. The setup of the three-point bending model is illustrated in Fig. 17. A notch with a depth of 70 mm and a width of $20\sqrt{3}$ mm is introduced at the center of the model. Both ends of the model are simply supported, and a prescribed displacement loading is applied to the central node on the opposite side of the notch. It is anticipated that strain localization will occur near the notch, while a more uniform displacement field will form in the regions far away from it.

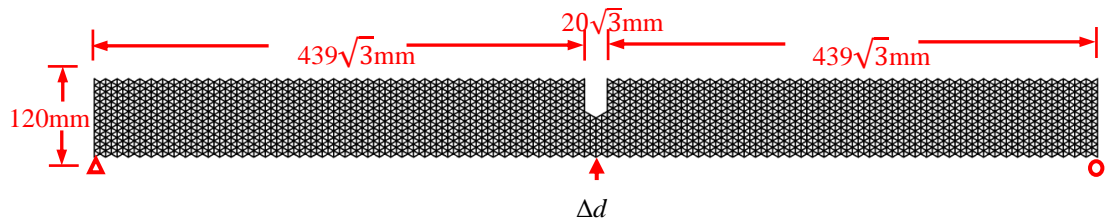

**Fig. 17.** Schematic illustration of the notched triangular lattice under the three-point bending test.

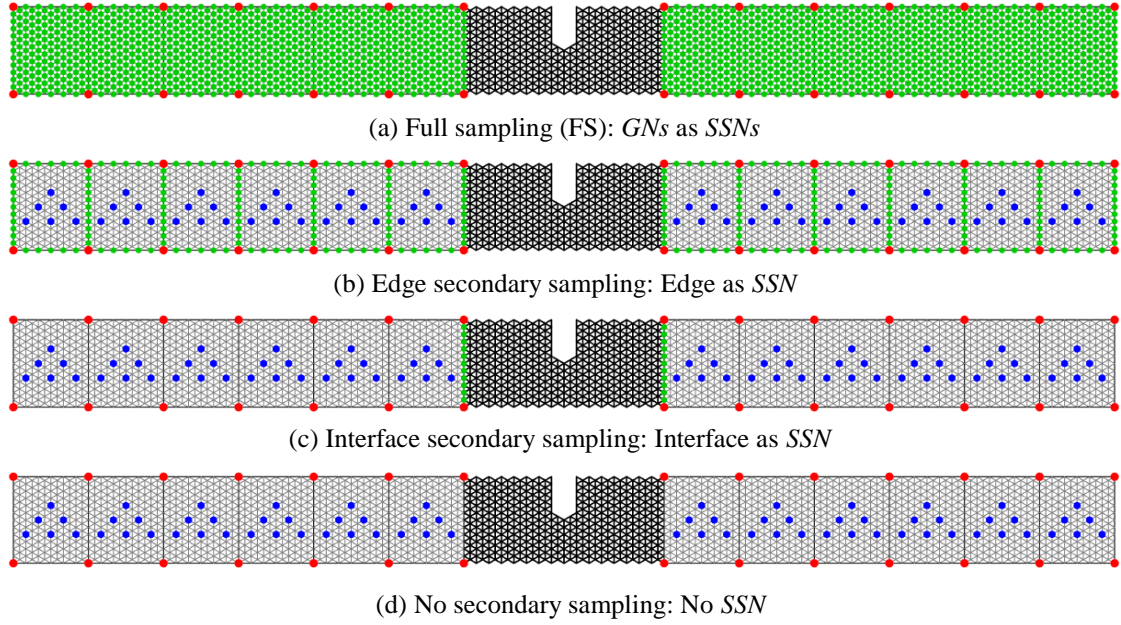

(a) Full sampling (FS): *GNs* as *SSNs*

(b) Edge secondary sampling: Edge as *SSN*

(c) Interface secondary sampling: Interface as *SSN*

(d) No secondary sampling: No *SSN*

**Fig. 18.** Schematic diagram of different sampling schemes for the three-point bending test. (a) Full sampling (FS): GNs as SSNs, (b) Edge secondary sampling: Edge as *SSN*, (c) Interface secondary sampling: Interface as *SSN* and (d) No secondary sampling: No *SSN*.

GNQC with four sampling schemes are applied to coarse-grain the notched lattice by employing bilinear element with a size of $12l_0$. GNQC models with full sampling (FS), edge secondary sampling (ESS), interface secondary sampling (ISS), and no secondary sampling (NSS), as shown in Fig. 18a-d, respectively. Fig. 19 plots the reaction force-displacement curves defined at the node where displacement loading is applied. At the end of the loading process, the reaction force gradually stabilizes, a behavior primarily caused by the progressive failure of truss elements near the notch and those surrounding the loaded nodes. Fig. 19 indicates a close match between full resolution (FR) results and the ones from the employed sampling schemes. Table 8 further quantifies the introduced error during the coarse-graining process for each sampling scheme by calculating the external work from the respective force-displacement curve using the trapezoidal method. For all the employed sampling schemes, the deviation of external work from the full resolution model is all kept within an acceptable range (less than 2.5%), showing the effectiveness of the proposed GNQC.

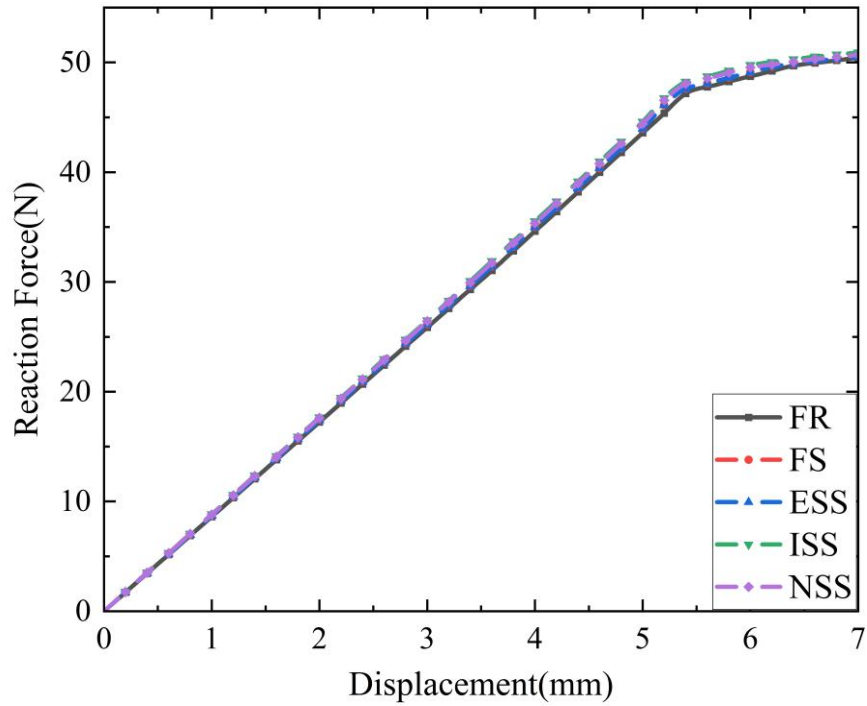

**Fig. 19.** Reaction force versus imposed mid-span displacement for the three-point bending. FR: full resolution model, FS: full sampling model, ESS: Edge secondary sampling, ISS: Interface secondary sampling and NSS: No secondary sampling.

**Table 8.** Comparison of full resolution model and models with different sampling schemes in terms of the numbers of sampling nodes and accuracy.

| Model | Number of sampling nodes | External work relative to FR | Reaction force at maximum strain relative to FR |
|---|---|---|---|
| FR | 2197 | 100% | 100% |
| FS | 2197 | 100.7% | 100.1% |
| ESS | 745 | 100.8% | 100.1% |
| ISS | 493 | 101.7% | 100.8% |
| NSS | 471 | 102.2% | 100.9% |

## 4.2 Uniaxial Tension

In this section, a uniaxial tensile simulation will be conducted on a triangular lattice model with a notch at the center along the *x*-axis (Fig. 20). The model consists of 3,253 lattice points with dimensions of $38l_0 \times 42\sqrt{3}l_0$, and displacement loadings are applied at the two ends (Fig. 20). GNQC is employed with the recommended interface secondary sampling (ISS) and bilinear elements of size $16l_0$ to coarse-grain the original model (Fig. 20b), as outlined in the previous study.

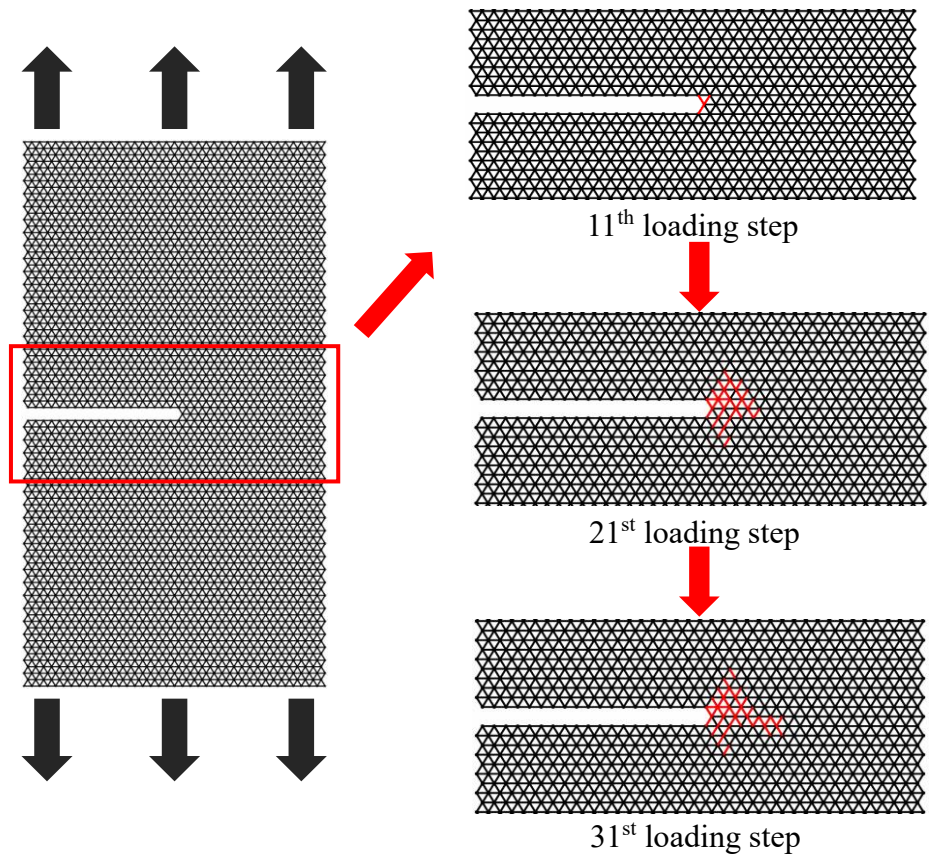

(a) Full resolution model

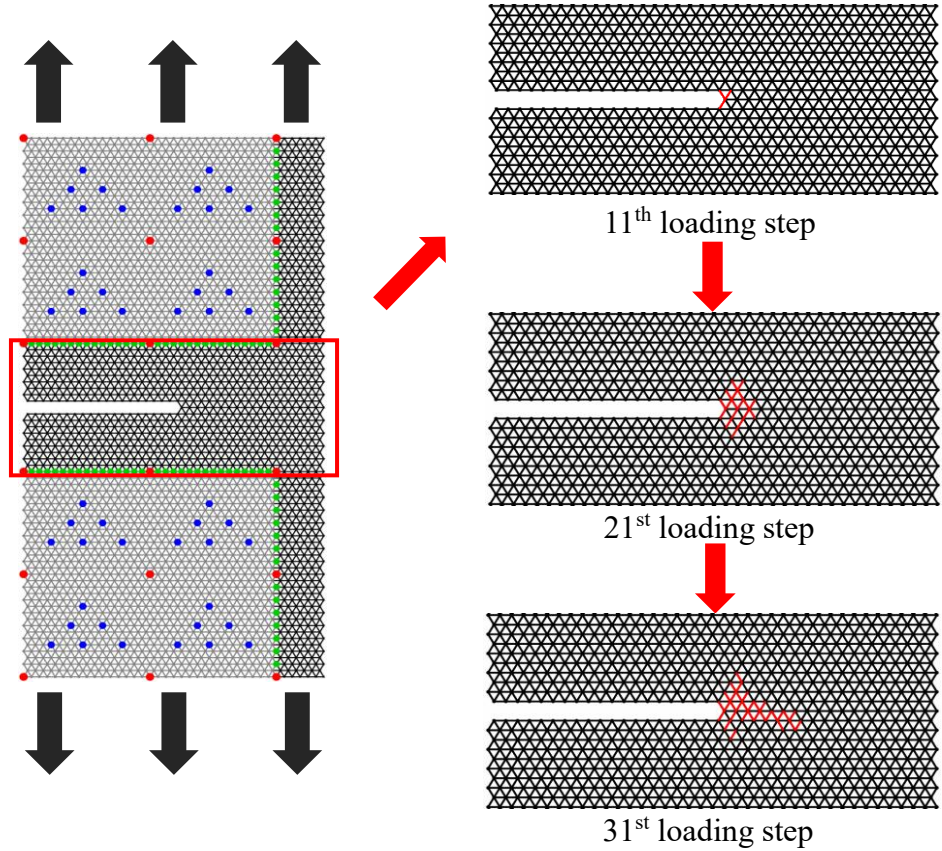

(b) GNQC model with interface secondary sampling

**Fig. 20.** Uniaxial tensile test of the notched lattice structure modeled by full resolution (FR) model (a) and GNQC with interface secondary sampling (b). The red connections between nodes represent the failed struts, while the black connections indicate intact lattice struts.

The local region near the notch is zoomed in and displayed in the right panels of Fig. 20, where the broken struts are represented in red. As shown, the crack initiates at the notch tip and gradually propagates into the surrounding region. GNQC with interface secondary sampling (ISS) accurately predicts the crack path, closely matching the one observed in the full resolution model at the corresponding loading steps. Furthermore, the force–displacement curve analysis presented in Fig. 21 quantitatively demonstrates the method's ability to capture crack propagation, further validating the effectiveness of the proposed GNQC in modeling the mechanical behavior of architected lattice structures.

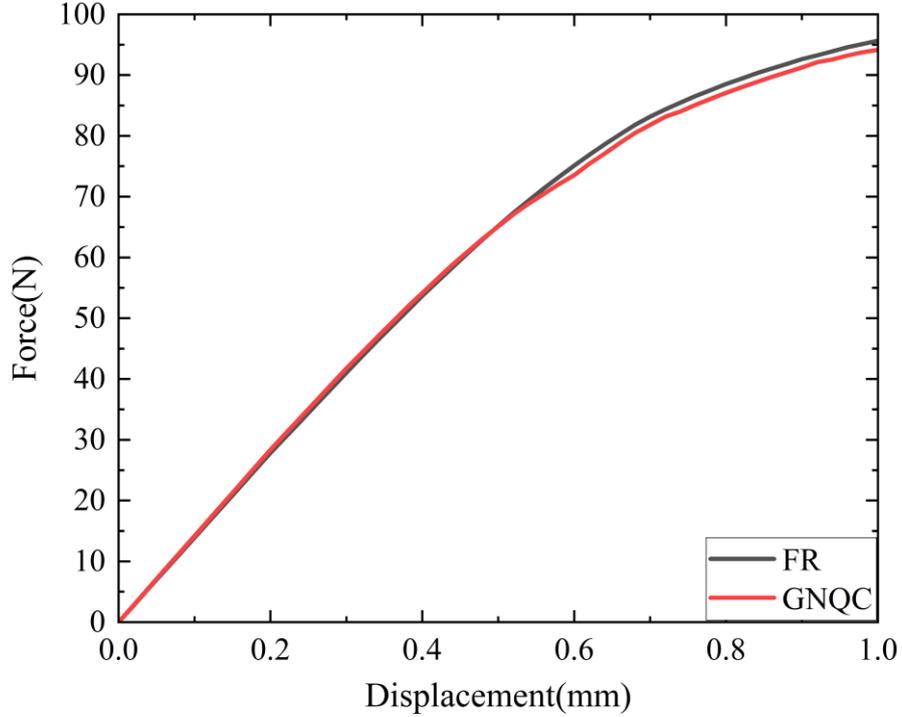

**Fig. 21.** The force–displacement curves from the respective full resolution (FR) model and GNQC with interface secondary sampling in crack propagation simulations.

# 5. Conclusions

We introduced the Generalized Non-local Quasicontinuum (GNQC) method as a concurrent multiscale approach for efficiently modeling the mechanical behavior of architected large-scale lattice structures. GNQC generalizes the classical nonlocal Quasicontinuum framework by eliminating the assumption of affine or high-order deformation patterns in energy sampling, ensuring consistency with general finite element shape functions. The method is characterized by the following three key features:

(1) Constitutive-Model-Consistent Framework: GNQC directly employs a discrete constitutive model that accurately captures the interactions between lattice nodes across the entire domain. This approach ensures monotonic convergence to the full-resolution model as the finite element size in the coarse-graining process approaches the lattice spacing.

(2) Shape-Function-Consistent Energy Sampling: The energy sampling strategy in GNQC is designed to align with the interpolation order of the employed finite element shape functions, thereby optimizing both accuracy and computational efficiency. This is achieved by mathematically determining the energy distribution pattern within the coarse-grained region, a method applicable to general finite element shape functions.

(3) Consistent Interfacial Compatibility: Seamless transfer of energy and force is maintained across the interface between the coarse-grained region and the locally full-resolution domain, eliminating the need for cumbersome interfacial treatments. This is facilitated by the use of complementary secondary sampling.

The performance of GNQC has been validated through multiple numerical tests—including tension, clamped bending, three-point bending, and crack propagation problems. The results closely match those from full-resolution models while significantly reducing computational costs. Future research will explore the application of GNQC with advanced beam theories as the constitutive model between lattice nodes, the development of adaptive parallel algorithms, and the extension of the method to dynamic scenarios such as impact and other dynamic loading conditions, as well as multi-lattice structures.

**Acknowledgements**

This work was supported by the following funding sources: the Natural Science Foundation of China (Grants 12272214 and 12472077), the Science and Technology Commission of Shanghai Municipality (Grant 23010500400), the Shanghai Gaofeng Project for University Academic Program Development, and the State Key Laboratory of Structural Analysis, Optimization and CAE Software for Industrial Equipment (Grant GZ23107).